\documentclass[12pt]{amsart}    
\usepackage{amsmath,amsthm,amssymb, amscd}

\newcommand{\uc}{{\mathcal U}}  
\newcommand{\acw}{\widetilde{\mathcal A}}      
\newcommand{\bc}{{\mathcal B}}    
\newcommand{\Ss}{{\mathcal S}}    
\newcommand{\ttheta}{\widetilde{\theta}}    
\newcommand{\kc}{{\mathcal K}}

\renewcommand{\ker}{\mbox{\rm Ker\,}}

\newcommand{\ku}{{\Bbbk}}    
\newcommand{\Z}{{\mathbb Z}}

\newcommand{\C}{{\mathbb C}}    
    
\newcommand{\VGamma}{\widehat{\Gamma}}  
\newcommand{\BGamma}{{\mathbb G}}
\newcommand{\bVGamma}{\widehat{{\BGamma}}}  
\newcommand{\ydg}{{}_{\Gamma}^{\Gamma}\mathcal{YD}}
\newcommand{\toba}{{\mathfrak B}}
\newcommand{\wtoba}{\widehat{\mathfrak B}}
\newcommand{\ftoba}{{\mathcal K}}

\numberwithin{equation}{section}

\newcommand{\End}{\mbox{\rm End\,}}

\newcommand{\id}{\mathop{\rm id\,}}     
    
\newcommand{\ord}{\mathop{\rm ord}}    
    
\newcommand{\ad}{\mbox{\rm ad\,}}    
    
\newcommand{\Alt}{\mbox{\bf S\,}}
\newcommand{\gr}{\mbox{\rm gr\,}}    
    
\newcommand{\abodrio}{{\mathcal A}\left(\Gamma,\, (a_{ij})_{1\le i, j\le \theta},\, 
(g_{i})_{1 \le i \le \theta}, \, (\chi_{j})_{1 \le j \le \theta},  
\, (\lambda_{ij})_{1 \le i < j \le \theta, \, i\not\sim j} \right)}

\theoremstyle{plain}    
    
\newtheorem{teo}[equation]{Theorem}    
\newtheorem{lema}[equation]{Lemma}    
\newtheorem{cor}[equation]{Corollary}    
\newtheorem{prop}[equation]{Proposition}

\theoremstyle{definition}    
\newtheorem{defi}[equation]{Definition}    
 \newtheorem{exa}[equation]{Example}    
    
\theoremstyle{remark}    
\newtheorem{obs}[equation]{Remark}
\newtheorem{rmk}[equation]{Remarks}

\def\pf{\begin{proof}}
\def\epf{\end{proof}}

\theoremstyle{remark}

\begin{document}

\renewcommand{\baselinestretch}{1.2}
\renewcommand{\thefootnote}{}
\thispagestyle{empty}
%\vspace*{2in}
\title{Finite quantum groups over abelian groups of prime exponent}
%Pointed Hopf algebras with abelian coradical of exponent $p$}
\author{ Nicol\'as Andruskiewitsch and Hans-J\"urgen Schneider}
\address{Facultad de Matem\'atica, Astronom\'\i a y F\'\i sica\\
Universidad Nacional de C\'ordoba \\ (5000) Ciudad Universitaria 
\\C\'ordoba \\Argentina}
\email{andrus@mate.uncor.edu}
\thanks{This paper is available at {\tt www.mate.uncor.edu/andrus/articulos.html}.
This work was partially supported by CONICET,
	CONICOR, DAAD, the Graduiertenkolleg of the Math. Institut 
	(Universit\"at \ M\"unchen)  and Secyt (UNC)}
\address{Mathematisches Institut\\
Universit\"at \ M\"unchen\\
Theresienstra\ss e 39\\
D-80333 M\"unchen\\
Germany}
\email{hanssch@rz.mathematik.uni-muenchen.de}
%\keywords{Pointed Hopf algebras, finite quantum groups}
\date{February 16. 2000}
\maketitle
    
\section{Introduction}

Since the discovery of quantum groups (Drinfeld, Jimbo) and  finite dimensional
variations thereof (Lusztig, Manin), these objects were studied from different points
of view and had many applications. The present paper is part of 
a series where we intend to show that  important classes of Hopf 
algebras are  quantum groups and therefore belong to Lie theory. 
One of our main results is the explicit 
construction of a general family of pointed Hopf algebras
from Dynkin diagrams (Theorem \ref{linking}). All the Frobenius-Lusztig kernels and their 
parabolic subalgebras belong to this family, but in addition we get many new examples. 
We show that any finite dimensional
pointed Hopf algebra with group of prime exponent (greater than 17) is indeed in
this family; 
see our main Theorem below. An important step in the proof follows
from a another main result (Theorem \ref{degree1}), where we show that a wide family of 
finite dimensional pointed Hopf algebras is generated by group-like 
and skew-primitive elements, giving additional support to a conjecture in \cite{AS2}.

\subsection{Main result\label{uno.1}}

We shall work over an algebraically closed field $\ku$ of characteristic 0. 
We denote by  $\widehat{\Gamma}$ the group of characters of an abelian group $\Gamma$. 
If $A$ is a Hopf algebra, then we denote by $G(A)$ the group of group-like
elements of $A$. 
Recall that $A$ is {\em pointed} if $\ku G(A)$ is the largest cosemisimple
subcoalgebra of $A$, or equivalently if
any irreducible $A$-comodule is one-dimensional. 

\bigbreak
Let $p$ be an odd prime number.  Let $s$ be a natural number and let
$\varGamma(s) = (\Z/(p))^s$. 

\begin{teo}\label{main} 
(a). Let $p> 17$.
Let $\mathcal A$ be a pointed finite-dimensional Hopf algebra such that
$G(\mathcal A) \simeq \varGamma(s)$. Then there exist 
\begin{itemize} \item a finite Cartan matrix 
$(a_{ij})\in \ku^{\theta\times \theta}$ \cite{K};
\item elements $g_{1}, \dots, g_{\theta} \in
\varGamma(s)$, $\chi_{1}, \dots, \chi_{\theta} \in
\widehat{\varGamma(s)}$ such that 
\begin{align}\label{diagnotriv}
\langle \chi_{i}, g_{i}\rangle &\neq 1, \qquad \text{for all }  1 \le i \le \theta, \\
\label{cartantype}
\langle \chi_{j}, g_{i}\rangle \langle \chi_{i}, g_{j}\rangle &=
\langle \chi_{i}, g_{i}\rangle^{a_{ij}}, \qquad \text{for all }  1 \le i, j \le \theta;
\end{align}
\item and a collection $(\lambda_{ij})_{1 \le i < j \le \theta, \, i\not\sim j}$ of elements 
in $\ku$ 
which is a "linking datum" for 
the previous data
as in Section \ref{link-dat} below;
\end{itemize}
such that $\mathcal A$ can be presented as algebra by generators $a_{1}, \dots, a_{\theta}$, 
$y_{1}, \dots, y_{s}$ and relations
\begin{align}\label{relations}
y_{h}^p &= 1, \qquad y_{m}y_{h} = y_{h}y_{m}, \qquad 1 \le m,h  \le s,\\\label{relations1}
y_{h}a_{j} &= \chi_{j}(y_{h})a_{j}y_{h}, \qquad 1 \le h \le s,\, 1 \le j \le \theta,\\
\label{relations2}
(\ad a_{i})^{1 - a_{ij}}a_{j} &= 0, \qquad 1 \le i \neq j \le \theta, \quad i \sim j, \\
\label{relations3}
a_{i}a_{j} - \chi_{j}(g_{i})a_{j}a_{i} &= \lambda_{ij}(1 - g_{i}g_{j}),  
\qquad 1 \le i < j \le \theta, \quad i \not\sim j;\\
\label{relations4}
a_{\alpha}^{p} &= 0,  
\qquad \alpha \in \Phi^{+};
\end{align}
and where the Hopf algebra structure is determined  by
\begin{equation}\label{sk-gl}
\Delta y_{h} = y_{h}\otimes y_{h}, \qquad
\Delta a_{i} = a_{i}\otimes 1 + g_{i} \otimes a_{i}, \qquad 1 \le h \le s,\, 1 \le i \le \theta.
\end{equation}

\bigbreak
(b). Conversely, let  
$(a_{ij})\in \Z^{\theta\times \theta}$
be a finite Cartan matrix, $g_{1}, \dots, g_{\theta} \in
\varGamma(s)$, $\chi_{1}, \dots, \chi_{\theta} \in
\widehat{\varGamma(s)}$ such that \eqref{diagnotriv}, \eqref{cartantype} hold
and $(\lambda_{ij})$  a "linking datum" for  
$(a_{ij})$, $g_{1}, \dots, g_{\theta}$ and $\chi_{1}, \dots, \chi_{\theta}$. 
Assume that $p>3$ if the Cartan matrix $(a_{ij})$ has a connected component
of type $G_2$.
Then the algebra $\mathcal A$ presented  
by generators $a_{1}, \dots, a_{\theta}$, 
$y_{1}, \dots, y_{s}$ and relations \eqref{relations}, \eqref{relations1}, 
\eqref{relations2},
\eqref{relations3},
\eqref{relations4}
has a unique Hopf algebra structure determined  by \eqref{sk-gl}. 
$\mathcal A$ is pointed, $G(\mathcal A) \simeq \varGamma(s)$ and $\dim \mathcal A 
= p^{s + \vert \Phi^{+}\vert}$.
\end{teo}

Here $\ad$ is the adjoint action; $i \sim j$, resp. $i \not\sim j$, means that $i$ and $j$ belong
to the same connected component, resp. to different connected components,
of the Dynkin diagram corresponding to $(a_{ij})$; 
$\Phi^{+}$ is the set of positive roots of the root system  associated to the Cartan matrix $(a_{ij})$; 
the "root vectors" $a_{\alpha}$ are defined in Subsection \ref{rootvectors}
below.

\begin{rmk}
(i). We define the notion of "linking datum" attached to a collection $g_{i}$,
$\chi_{i}$, $1\le i \le \theta$, $(a_{ij})_{1\le i, j \le \theta}$ 
for a general finite abelian group $\Gamma$, cf. Section  \ref{link-dat}.
It is always possible to consider
 linking data with entries $\lambda_{ij}$ taking
only the values 1 or 0. It follows then from Theorem \ref{main} that there are only
finitely many isomorphism classes of
finite dimensional Hopf algebras with coradical $\ku
\varGamma(s)$. For more general finite abelian groups, this is no longer
true \cite{AS1}, \cite{BDG}, \cite{G}. 

\bigbreak
(ii). Let $(a_{ij})_{1\le i, j \le \theta}$ be a finite Cartan matrix. The
problem of finding all the collections $g_{i}\in \varGamma(s)$, $\chi_{i}\in
\widehat{\varGamma(s)}$, $1\le i \le \theta$, such that \eqref{diagnotriv}
and \eqref{cartantype}
hold has been discussed in \cite{AS2}.  It can be stated as the problem of
finding all the solutions of a system of algebraic equations over $\Z /(p)$
and it is in principle solvable. Note that in particular
$$\theta \le 2s\frac{p-1}{p-2},
$$
see \cite[Prop. 8.3]{AS2}.
Next, the question of finding all the possible
linking data  attached to a fixed collection $g_{i}$,
$\chi_{i}$, $1\le i \le \theta$, $(a_{ij})_{1\le i, j \le \theta}$, 
has also a strong combinatorial flavor, see  Section  
\ref{link-dat}.
Once these two problems are solved effectively, the determination of all
isomorphism classes is an accessible task using \cite[Prop. 6.3]{AS2}, \cite[Lemma 1.2]{AS3}.

\bigbreak
(iii). In particular, the complete classification of all finite
dimensional pointed Hopf algebras with group of group-likes
$\varGamma(1) = \Z/(p)$, $p\neq 5, 7$ follows from Theorem \ref{main}: it is the list given in
\cite[Theorem 1.3]{AS2} plus the Frobenius-Lusztig kernels as described in \cite{AS1}.
Indeed, replacing in the proof of Theorem \ref{main}  \cite[Cor. 1.2]{AS2}
by \cite[Th. 1.3]{AS2} we get the classification for all primes $p \neq 5$ or  $7$, 
in view of Theorem \ref{QSR} and \cite[Lemma 4.2]{AS3}. The only cases not covered
are $p = 5$, type $B_2$ and  $p=7$, type $G_2$.
This result was independently obtained by Musson \cite{Ms}
using different methods starting from our previous article \cite{AS2}.
%; he  obtained the classification also for $p = 5, 7$
%showing in particular that the quantum Serre relations do not hold for these 
%small primes.

\bigbreak
(iv). Part (b) of Theorem \ref{main} follows
from a more general result for arbitrary finite abelian groups, see Theorem \ref{linking}.
For suitable choices of $(a_{ij})_{1\le i, j \le \theta}$, $g_{i}$
and $\chi_{i}$, $1\le i \le \theta$, and $(\lambda_{ij})_{1\le i, j \le \theta, \, i\not\sim j}$ 
one recovers Frobenius-Lusztig kernels
and their parabolic subalgebras. See example \ref{FLK}.
Otherwise  Theorem \ref{linking} provides
many new examples of finite dimensional Hopf algebras arising from exotic linking data.

\bigbreak
(v). The dimensions of the Hopf algebras in Theorem \ref{uno.1} are very special numbers. This phenomenon is shown in general for arbitrary finite groups $\Gamma$ in Theorem \ref{dimension}.

\bigbreak
(vi). Up to now,  the determination of all finite dimensional pointed Hopf
algebras $A$ with $G(A) \simeq \Gamma$,  for a fixed group $\Gamma$, was
known only for $\Gamma = \Z/(2)$ \cite{N}. The case of $\Gamma$ of
exponent 2 will be treated in a separate article with S. D\u{a}sc\u{a}lescu.  
Other classification results
of pointed Hopf algebras are known for some fixed dimension $d$: $d = p^{2}$
is easy and follows from \cite{N}, \cite{NZ}; $d = p^{3}$ was done in
\cite{AS1},  and by different methods in \cite{CD}, \cite{SvO};  $d = p^{4}$
in \cite{AS3} (and does not seem to be possible via the other
methods); $d = 16$ in \cite{CDR}, $d = 32$ in \cite{Gn1}.

\bigbreak
(vii). The classification  of all {\it coradically graded}
pointed Hopf algebras of dimension 
$p^{5}$ was obtained in \cite{Gn2}. It is not difficult to deduce the 
classification  of all pointed Hopf algebras of dimension 
$p^{5}$ using Theorem \ref{main}
and results in \cite{AS3}.

\bigbreak
(viii). The results of this paper heavily depend on our  paper \cite{AS2}
and on previous work on quantum groups \cite{L1}, \cite{L2}, \cite{Ro1},
\cite{Ro2}, \cite{dCP}, \cite{Mu}.

\end{rmk}

\subsection{Conventions\label{dos.1}}
 Our reference for the theory of Hopf algebras is 
\cite{Mo}.
The notation for Hopf algebras is  standard:
$\Delta$, $\Ss$, $\epsilon$, denote respectively the
comultiplication, the antipode, the counit; we use Sweedler's
notation but dropping the summation symbol.   

If $C$ is a coalgebra then $G(C)$ denotes the set of its group-like elements 
and $C_{0} \subset C_{1} \subset \dots$ its coradical filtration. So that
$C_{0}$ is the coradical of $C$. As usual, ${\mathcal
P}_{g,h}(C)$ denotes the space of $(g,h)$-skew primitives of $C$, $g,h 
\in G(C)$. If $C$ has a distinguished group-like 1, then
we set ${\mathcal P}(C) = {\mathcal P}_{1,1}(C)$, 
the space of primitive elements of $C$.

If $A$ is an algebra and $(x_{i})_{i\in I}$ is a family of 
elements of $A$ then $\ku\langle x_{i}\rangle_{i\in I}$ or simply
$\ku\langle x_{i}\rangle$, resp. $\langle x_{i}\rangle_{i\in I}$ or 
$\langle x_{i}\rangle$
denotes the subalgebra, resp. the two-sided ideal, generated by the
$x_{i}$'s.

Let  $H$ be  a Hopf algebra. 
A Yetter-Drinfeld  module over $H$ 
is a vector space $V$ provided 
with structures of left $H$-module and left $H$-comodule such that 
$\delta(h.v) = h_{(1)} v_{(-1)} \mathcal S h_{(3)} \otimes h_{(2)}.v_{(0)}.$
We denote  by $_{H}^{H}\mathcal{YD}$
the (braided) category of Yetter-Drinfeld  modules over $H$.

Assume that  $H = \ku\Gamma$ where $\Gamma$ is a finite abelian group.
We denote $_{\Gamma}^{\Gamma}\mathcal{YD} := _{H}^{H}\mathcal{YD}$.
Let $g \in \Gamma$, $\chi \in \VGamma$ and $V$ a module, resp. a comodule,
resp. a Yetter-Drinfeld module,
over $\Gamma$. Then we denote $V^{\chi} 
= \{v \in V: h.v = \chi(h)v, \quad \forall h\in \Gamma\}$, resp.  
$V_{g} = \{v \in V: \delta(v) = g\otimes v\}$,  
resp. $V_{g}^{\chi} := V_{g}\cap V^{\chi}$.  
If $V$ is a locally finite Yetter-Drinfeld module,
then $V = \oplus_{g \in \Gamma, \chi \in \VGamma} V_{g}^{\chi}$.
Conversely, a vector space $V$ provided with a direct sum decomposition 
$V = \oplus_{g \in \Gamma, \chi \in \VGamma} V_{g}^{\chi}$
has an evident Yetter-Drinfeld module structure.

\section{Outline of the paper and proof of the main result}
Theorem \ref{main}
follows from Theorems \ref{flk-tw}, \ref{linking}, \ref{QSR}, \ref{genlifting} and Corollary \ref{degree1p}
in the present article, along the guidelines proposed in \cite{AS1}. 
We give now the proof of Theorem \ref{main}
assuming those results which hold over arbitrary finite abelian groups. 
This section serves also as a guide to the different sections of  
the paper.

\subsection{The proof}
Let $\mathcal A$ be a finite dimensional pointed Hopf algebra with $G(\mathcal A) \simeq
\varGamma(s)$. Let 
$\gr \mathcal A := \oplus_{n \ge 0} \gr \mathcal A(n)$,
where $\gr A(0) = A_{0}$, $\gr \mathcal A(n) = \mathcal A_{n}/\mathcal A_{n-1}$, if $n>0$
be the graded coalgebra associated to the coradical filtration of $\mathcal A$. 
Then $\gr \mathcal A$ is a graded Hopf algebra \cite{Mo} and both the inclusion
$\iota: \mathcal A_{0} \hookrightarrow \gr \mathcal A$ and the graded projection $\pi: \gr 
\mathcal A
\to \mathcal A_{0}$ are Hopf algebra maps. Let $R := \gr \mathcal A^{co\pi} = \{x\in \gr 
\mathcal A:
(\id\otimes \pi)\Delta(x) = x\otimes 1\}$; it is a graded braided Hopf
algebra in $^{\varGamma(s)}_{\varGamma(s)}\mathcal {YD}$ with the grading
inherited from $\gr \mathcal A$: $R = \oplus_{n \ge 0} R(n)$, $R(n) := R\cap \gr \mathcal A(n)$.
Notice that $\gr \mathcal A$ can be reconstructed from $R$ as a biproduct:
$$\gr \mathcal A \simeq R \# \ku \varGamma(s).$$

The braided Hopf algebra $R$ is called the {\em diagram} of $\mathcal A$. One has  
\begin{flalign}\label{nicholsuno} &R(0) = \ku 1, &\\
\label{nicholsdos} &R(1) = {\mathcal P}(R),&
\end{flalign}
and we know from Corollary \ref{degree1p} below that
\begin{flalign}\label{nicholstres} &R \text{ is generated as an algebra by R(1)}.&  \end{flalign}
Let $V := R(1)$; it is a Yetter-Drinfeld submodule of $R$. Since $R$
satisfies \eqref{nicholsuno}, \eqref{nicholsdos} and \eqref{nicholstres} 
we know that $R \simeq \toba(V)$ is a Nichols
algebra, see Subsection \ref{dos.nichalg}.
Now there exists a basis $x_{1}, \dots, x_{\theta}$ of $V$
and $g_{1}, \dots, g_{\theta} \in \varGamma(s)$, $\chi_{1}, \dots,
\chi_{\theta}\in \widehat{\varGamma(s)}$ such that $x_{i} \in V_{g_{i}}^{\chi_{i}}$,
$1\le i \le \theta$. 
Since $\mathcal A$ is finite dimensional,
$\langle \chi_{i}, g_{i}\rangle \neq 1$ for all $i$  \cite[Lemma
3.1]{AS1}
and there is a finite Cartan matrix $(a_{ij})_{1\le i, j\le \theta}$ such that
\eqref{cartantype}
holds  \cite[Cor. 1.2]{AS2}.

\bigbreak To give an explicit description of $\toba(V)$, 
we introduce {\em root vectors} in $\toba(V)$ generalizing the
root vectors defined in \cite{L1}. We note that Lusztig's root vectors 
can be described up to a non-zero scalar as an iterated braided commutator of 
simple root vectors. We then define the root vectors in the general case by
exactly the same iterated braided commutator with respect to our more general
braiding.  
As one of our main results, 
we obtain a presentation by generators and relations and a PBW
basis for $\toba(V)$ from the corresponding Theorem for Frobenius-Lusztig
kernels, using Drinfeld's twisting essentially in the same way as in
\cite{AS2}. See Theorem \ref{flk-tw}.  We can then deduce part (b)
of Theorem \ref{main}: for connected Dynkin diagrams it
is a consequence of Theorem \ref{flk-tw}; the non-connected case is dealt with
twisting of the algebra structure \cite{DT}. See Theorem \ref{linking}.

It follows at once from Theorem \ref{flk-tw}
that $\gr \mathcal A$ can be presented as algebra by generators 
$y_{1}, \dots, y_{s}$ (homogeneous of degree 0) and 
$x_{1}, \dots, x_{\theta}$ (homogeneous of degree 1), 
and relations
\begin{align}\label{gr-relations}
y_{h}^p &= 1, \qquad y_{m}y_{h} = y_{h}y_{m}, \qquad 1 \le m,h  \le s,\\\label{gr-relations1}
y_{h}x_{j} &= \chi_{j}(y_{h})x_{j}y_{h}, \qquad 1 \le h \le s,\, 1 \le j \le \theta,\\
\label{gr-relations2}
(\ad x_{i})^{1 - a_{ij}}x_{j} &= 0, \qquad 1 \le i \neq j \le \theta;\\
\label{gr-relations4}
x_{\alpha}^{p} &= 0,  
\qquad \alpha \in \Phi^{+};
\end{align}
and where the Hopf algebra structure is determined  by
\begin{equation}\label{gr-sk-gl}
\Delta y_{h} = y_{h}\otimes y_{h}, \qquad
\Delta x_{i} = x_{i}\otimes 1 + g_{i} \otimes x_{i}, \qquad 1 \le h \le s,\, 1 \le i \le \theta.
\end{equation}

By \cite[Lemma 5.4]{AS1}, we can choose $a_{i} \in {\mathcal P}_{g_{i}, 1}(\mathcal A)^{\chi_{i}}$ 
such that the class of $a_{i}$ in $\gr \mathcal A(1)$ coincides with $x_{i}$. We also keep the notation
$y_{j}$ for the generators of $G(\mathcal A)$. It is clear that relations
\eqref{relations} and \eqref{relations1} hold. Now relations
\eqref{relations2} and \eqref{relations3}, resp. \eqref{relations4}, 
hold because of Theorem \ref{QSR},
resp. Lemma \ref{powrootvec}. 

The Theorem now follows from Theorem \ref{genlifting}. \qed

\subsection{What is next?}
\noindent There are several obstructions to extend Theorem \ref{main}
to general finite abelian groups. First, it
is  open whether the diagram of a finite dimensional pointed Hopf algebra 
is generated in degree one, {\em i.e.} when it is a
Nichols algebra; second, there are finite dimensional Nichols 
algebras which are not of Cartan type \cite{N}. 

For liftings of $\gr A$ when $R$ is a Nichols algebra of Cartan type,
quantum Serre relations of connected vertices still hold as we show in
Theorem \ref{QSR}
below; however the powers of root vectors are not necessarily 0. We should
have  $a_{\alpha}^{N_{\alpha}} = u_{\alpha} \in \ku G(A)$; the determination
of $u_{\alpha}$ when $\alpha$ is a non-simple root will be addressed in the
next paper (it was done in \cite{AS3} for type $A_{2}$). 

\bigbreak
\section{Braided Hopf algebras \label{dos}}

\subsection{Biproducts\label{dos.2}}
Let $R$ be a braided Hopf algebra in $_{H}^{H}\mathcal{YD}$;
this means that $R$ is an algebra and a coalgebra in $_{H}^{H}\mathcal{YD}$
and that the comultiplication $\Delta_{R}: R \to R\otimes R$
is an algebra map when in $R\otimes R$ the multiplication {\em twisted
by the braiding } $c$ is considered; in addition $R$ admits an antipode.
To avoid confusions we use the following variant of Sweedler's notation 
for the comultiplication of $R$:
$\Delta_R(r) = r^{(1)} \otimes r^{(2)}$.
Let $A = R\# H$ be the biproduct or bosonization of $R$ \cite{Mj}, \cite{Ra}. 
Recall that the multiplication and comultiplication of $A$ are given by
$$
(r\# h)(s\# f) = r (h_{(1)}.s)\# h_{(2)}f, \quad \Delta(r\# h) = 
r^{(1)} \# (r^{(2)})_{(-1)} h_{(1)} \otimes (r^{(2)})_{(0)}
 \#  h_{(2)}. $$

The maps 
$\pi: A \to H$ and $\iota: H \to A$, $\pi(r\# h) = \epsilon(r)h$,   
$\iota(h) = 1\# h$, are Hopf algebra homomorphisms; we have
$R = \{a\in A: (\id\otimes \pi)\Delta (a) = a\otimes 1\}$. 
Conversely, let $A$, $H$ be Hopf algebras provided 
with Hopf algebra homomorphisms $\pi: A \to H$ and $\iota: H \to A$.
Then $R = \{a\in A: (\id\otimes \pi)\Delta (a) = a\otimes 1\}$
is a braided Hopf algebra in $_{H}^{H}\mathcal{YD}$. 
The action  $.$ of $H$ on $R$ 
is the restriction of the adjoint action
(composed with $\iota$) and the coaction is $(\pi \otimes \id)\Delta$;
$R$ is a subalgebra of $A$ and the comultiplication is $\Delta_R(r) =
r_{(1)} \iota\pi\mathcal S(r_{(2)}) \otimes r_{(3)}$. These constructions are 
inverse to each other. We shall mostly omit $\iota$ in what follows.

Let $\vartheta: A \to R$ be the map given by
$\vartheta(a) = a_{(1)} \pi\Ss(a_{(2)})$.
Then 
\begin{equation}\label{varthetadjoint}
\vartheta(ab) = a_{(1)}\vartheta(b) \pi\Ss(a_{(2)}),\end{equation} for all 
$a,b \in A$ and $\vartheta(h) = \varepsilon(h)$ for all $h\in H$; therefore,
for all $a\in A$, $h\in H$, we have
$\vartheta(ah) = \vartheta(a) \varepsilon(h)$ and
\begin{equation}\label{varthetamodule}
\vartheta(ha) = h.\vartheta(a) = \vartheta(h_{(1)} a \pi\Ss(h_{(2)})).
\end{equation}
Notice also that $\vartheta$ induces a coalgebra isomorphism $A/AH^{+} 
\simeq R$. In fact, the isomorphism $A \to R\# H$ can be expressed explicitly as
$$
a \mapsto \vartheta(a_{(1)}) \# \pi(a_{(2)}), \qquad a\in A.
$$

If $A$ is a Hopf algebra, the well-known adjoint representation $\ad$
of $A$ on itself is given by
$\ad x(y)= x_{(1)} y \Ss(x_{(2)})$.
If $R$ is a braided Hopf algebra in $_{H}^{H}\mathcal{YD}$ 
then there is also a braided adjoint 
representation $\ad_c$ of $R$ on itself given by
$$\ad_c x(y)=\mu(\mu\otimes\Ss)(\id\otimes c)(\Delta\otimes\id)(x\otimes y),$$
where $\mu$ is the multiplication and
$c\in\End(R\otimes R)$ is the braiding. 
Note that if $x\in {\mathcal P} (R)$ then the braided adjoint representation 
of $x$ is just 
\begin{equation}\label{br-adj}
\ad_cx(y)=\mu(\id-c)(x\otimes y) =: [x, y]_{c}.
\end{equation}
The element $[x, y]_{c}$ defined by the second equality for any $x$ and $y$, regardless
of whether $x$ is primitive, will be called a braided commutator.

When $A = R\# H$, then for all $b,d\in R$, 
\begin{equation}\label{bradj}
\ad_{(b\# 1)}(d\# 1)=(\ad_cb(d))\# 1.\end{equation}

\subsection{Nichols algebras\label{dos.nichalg}}    
Let  $H$ be  a Hopf algebra and let $R = \oplus_{n \in {\mathbb N}} R(n)$ 
be a {\em graded} braided Hopf algebra in $_{H}^{H}\mathcal{YD}$. 
We say that $R$ is a {\em Nichols algebra} if
\ref{nicholsuno}, \ref{nicholsdos} and \ref{nicholstres} hold,
{\em cf.} \cite{N}, \cite{AS2}, \cite{AG}.
A Nichols algebra $R$ is uniquely determined by the Yetter-Drinfeld
module ${\mathcal P}(R)$; given a Yetter-Drinfeld module $V$, there exists a unique 
(up to isomorphism)  Nichols algebra $R$ with ${\mathcal P}(R) \simeq V$. 
It will be denoted $\mathfrak B(V)$. In fact, the kernel of
the canonical map $\varpi: T(V) \to \mathfrak B(V)$ can be described 
in several different ways. For instance, $\ker \varpi = \oplus_{n \geq 0}
\ker \Alt_{n}$ where $\Alt_{n}$ is the "quantum antisymmetrizer"
defined from the braiding $c$;
so that $\mathfrak B(V)$ is a "quantum shuffle algebra" 
and as algebra and coalgebra
only depends on the braiding $c: V\otimes V \to V\otimes V$.
See \cite{N}, \cite{Wo}, \cite{L3}, \cite{Ro1}, \cite{Ro2}, \cite{Sch}.

 Let $H = \ku\Gamma$ where $\Gamma$ is a finite abelian group. 
 Let $V$ be a finite dimensional Yetter-Drinfeld module over $\Gamma$.
Then there exist a basis $x_{1}, \dots, x_{\theta}$
of $V$ and elements $g_{1}, \dots, g_{\theta} \in \Gamma$, 
$\chi_{1}, \dots, \chi_{\theta} \in \VGamma$ such that
\begin{equation}\label{base}
x_{j} \in V_{g_{j}}^{\chi_{j}}, \qquad 1\le j \le \theta.
\end{equation}
In what follows we shall only consider
Yetter-Drinfeld modules $V$ such that 
$\langle \chi(i), g(i)\rangle \neq 1$, $1\le i \le \theta$.
The braiding $c$ is given with respect to the basis $x_{i}\otimes x_{j}$  by 
$c(x_{i}\otimes x_{j}) = b_{ij} \, x_{j}\otimes x_{i}$, 
where $$(b_{ij})_{1 \le i, j\le \theta}
=(\langle \chi(j), g(i)\rangle)_{1 \le i, j\le \theta}.$$

\begin{obs}\label{nodep}
Let $V$, resp. $\widetilde V$, be a finite dimensional Yetter-Drinfeld module over $\Gamma$, resp. $\widetilde \Gamma$,  with a basis
$x_{1}, \dots, x_{\theta}$ such that $x_i\in V^{\chi_i}_{g_i}$, resp. with a basis
$\widetilde x_{1}, \dots, \widetilde x_{\theta}$ such that $\widetilde x_i\in V^{\widetilde \chi_i}_{\widetilde g_i}$. Assume that $\chi_i(g_j) = \widetilde \chi_i(\widetilde g_j)$ for all $1\le i,j\le \theta$. Then there exists a unique algebra and coalgebra isomorphism ${\mathfrak B}(V)
\to {\mathfrak B}(\widetilde V)$ such that $x_i \mapsto \widetilde x_i$ for all
$1\le i\le \theta$.
\end{obs}

\begin{defi}\label{base-cartan} We shall say that a braiding given by a matrix 
$\mathbf b = (b_{ij})_{1\leq i,j \leq \theta}$  whose entries are roots of unity
  is {\it of Cartan type} 
if for all  $i, j$, $b_{ii}\neq 1$ 
  and   there exists  $a_{ij} \in \Z$ such that 
$$ b_{ij}b_{ji} = b_{ii}^{a_{ij}}.  $$ 
The integers $a_{ij}$ are uniquely determined by the following rules: 
\begin{itemize} \item If $i = j$ we take $a_{ii} = 2$; 
\item     
 if $i\ne j$, we select the unique $a_{ij}$ such that 
$-\ord b_{ii} < a_{ij} \leq 0$. \end{itemize}
\end{defi}

Then $(a_{ij})$ is a generalized Cartan matrix \cite{K}.
We shall say a  Yetter-Drinfeld module $V$ is of Cartan type,
resp. finite Cartan type, if its corresponding braiding
is of Cartan type, resp. the same plus the matrix
$(a_{ij})$ is of finite type.

\subsection{The twisting functor\label{dos.3}}
Let $H$ be a Hopf algebra and let $F$ be an invertible element
in $H\otimes H$ such that
\begin{equation}\label{twist}
F_{12} (\Delta \otimes \id) F = F_{23} (\id \otimes \Delta) F, \qquad
(\varepsilon \otimes \id)  (F) = 1 = (\id \otimes \varepsilon)  (F).
\end{equation} 
Then $H_{F}$, the same algebra $H$ 
but with the comultiplication $\Delta_{F} := F\Delta F^{-1}$, is again a Hopf algebra \cite{Dr}.    We shall write $F = F^{1}\otimes F^{2}$,
$F^{-1} = G^{1}\otimes G^{2}$; the new comultiplication will be denoted by
$\Delta_F(h) = h_{(1, F)} \otimes h_{(2, F)}$.

Let now $R$ be a braided Hopf algebra in $_{H}^{H}\mathcal{YD}$,
let $A = R\# H$ be its bosonization and 
consider the Hopf algebra $A_F$.
It follows from the definitions that $\pi: A_{F} \to H_{F}$ and 
$\iota: H_{F} \to A_{F}$ are also Hopf algebra homomorphisms. Hence 
$$R_{F} := \{a\in A_{F}: (\id\otimes \pi)\Delta_{F} (a) = a\otimes 1\}$$ 
is a braided Hopf algebra in the category $_{H_{F}}^{H_{F}}\mathcal{YD}$.  
We consider the corresponding map $\vartheta_{F}$ and define
$\psi: R \to R_{F}$ by
\begin{equation} \label{defpsi}
 \psi(r) = \vartheta_{F}(r), \qquad r\in R. \end{equation} 
The map $\psi$ was defined in \cite{AS2} in the case $H = \ku\Gamma$ is the group
algebra of a finite abelian group. The following Lemma generalizes
\cite[Lemma 2.3]{AS2}; part (iii), new even for $H = \ku\Gamma$,
will be needed in the sequel.

\begin{lema}\label{propspsi} (i). $\psi$  is an isomorphism of  $H$-modules. (Recall
that $H = H_{F}$ as algebras). 
%The behaviour of $\psi$ with respect to the comodule structures is ??

\medbreak
(ii). If $r, s \in R$ then  
\begin{equation} \label{2.7}
 \psi(rs) = F^{1}.\psi(r)\, F^{2}.\psi(s). \end{equation}  

\medbreak (iii).  If $r \in R$ then  
\begin{equation} \label{2.8}
 \Delta_{R_F}\psi(r) = F^{1}.\psi(r^{(1)})\otimes F^{2}.\psi(r^{(2)}). \end{equation}

%(iv). If $r, s \in R$ then  
%\begin{equation} \label{2.9}
%(\psi \otimes \psi) c(r \otimes s) = F. c_{F}(\psi(r) \otimes\psi(s)). 
%\end{equation}

\medbreak (iv). If $R$ is a graded braided Hopf algebra, then $R_{F}$ also is
 and $\psi$ is a graded map. If $R$ is a coradically graded braided Hopf 
algebra (resp. a Nichols algebra), then $R_{F}$ also is. \end{lema}  

\pf (i) follows from (\ref{varthetamodule}): 
$
\psi(h.r) = \vartheta_F(h.r) =  \vartheta_F(hr) =  h.\vartheta_F(r) = 
h.\psi(r).$
 Now we prove (ii):
\begin{align*}
\psi(rs) &=\vartheta_F(rs) = r_{(1, F)} \vartheta_F(s) \pi(\Ss_{F}(r_{(2, F)})) \\
&= r_{(1, F)} \pi(\Ss_{F}(r_{(2, F)})) \pi(r_{(3, F)})\vartheta_F(s) 
\pi(\Ss_{F}(r_{(4, F)})) \\ &= \psi(r_{(1, F)}) \pi(r_{(2, F)}). \psi(s) =
\psi(F^{1}r_{(1)}G^{1}) \pi(F^{2}r_{(2)}G^{2}). \psi(s) 
\\ &= \psi(F^{1}r_{(1)}) \varepsilon(G^{1}) \pi(F^{2}r_{(2)}G^{2}). \psi(s)
= F^{1}.\psi(r_{(1)}) \pi(F^{2})\pi(r_{(2)}). \psi(s)\\
&= F^{1}.\psi(r) \pi(F^{2}). \psi(s),
\end{align*}
as claimed. Here we have used (\ref{varthetadjoint}), the definitions
and (\ref{twist}). For the proof of (iii), we first observe that, if $r\in R$,
then
\begin{equation}\label{varthetacoalg}
\psi(r^{(1)}) \otimes \psi(r^{(2)}) =\vartheta_{F}(r_{(1)} 
\pi \mathcal S(r_{(2)}))
\otimes \vartheta_{F}(r_{(3)}) = \psi(r_{(1)}) \otimes \psi(r_{(2)}).
\end{equation}
Using  that $\vartheta_{F}$ is a coalgebra map,
(\ref{varthetamodule}) and (\ref{varthetacoalg}), we conclude that

\begin{align*}
 \Delta_{R_F}\psi(r) &= \Delta_{R_F}\vartheta_F(r) = \vartheta_F(r_{(1, F)})
\otimes \vartheta_F(r_{(2, F)})
 \\  &=  \vartheta_F(F^{1}r_{(1)}G^{1})
\otimes \vartheta_F(F^{2}r_{(2)}G^{2}) = \vartheta_F(F^{1}r_{(1)})
\otimes \vartheta_F(F^{2}r_{(2)}) 
\\ &= F^{1}.\vartheta_F(r_{(1)})
\otimes F^{2}.\vartheta_F(r_{(2)}) = F^{1}.\psi(r^{(1)})
\otimes F^{2}. \psi(r^{(2)}) ;
\end{align*}

The proof of (iv) has no difference with the proof of the analogous statement in
\cite[Lemma 2.3]{AS2}.
\epf

\bigbreak We now consider the special case when 
$H = \ku \Gamma$, $\Gamma$ a finite abelian group.  
Let $\omega: \VGamma\times \VGamma \to \ku^{\times}$ be a 2-cocycle,  i.e. 
$\omega(\tau, 1) = \omega(1, \tau) = 1$ and  $\omega(\tau, \zeta)\omega(\tau\zeta, \eta) = \omega(\tau, \zeta\eta) \omega(\zeta, \eta)$.
The cocycle $\omega$ allows to define a map 
$\Psi: \VGamma \times \Gamma \to \Gamma$ by
\begin{equation} \label{2.81}
 \langle \tau, \Psi(\chi, g)\rangle =  
\omega(\tau, \chi)\omega(\chi,\tau)^{-1}\langle \tau,  g\rangle, 
\qquad \tau \in \VGamma.
\end{equation}  
We identify $H$ with  the Hopf algebra $\ku^{\VGamma}$ of functions on the 
group $\VGamma$; we denote by $\delta_{\tau}\in H$ the function given by $\delta_{\tau}(\zeta) 
= \delta_{\tau, \zeta}$, $\tau, \zeta\in \VGamma$. Then $\delta_{\tau} = 
\frac 1{\vert \Gamma \vert}\sum_{g\in \Gamma} \langle \tau,  g^{-1}\rangle g.$  
Let $F \in H\otimes H$ be given by  
$$F = \sum_{\tau, \zeta\in \VGamma} \omega(\tau, \zeta) \delta_{\tau} \otimes \delta_{\zeta}. $$
Then $F$ satisfies \eqref{twist}; note that $H = H_{F}$. 
Let now $R$ be a braided Hopf algebra
in $\ydg$; we can consider the Hopf algebras $A = R\# \ku \Gamma$
and $A_{F}$, the braided Hopf algebra $R_{F} \in \ydg$ and
the map $\psi: R \to R_{F}$. We have
$$\psi (r) = \sum_{\tau\in \VGamma} \omega(\chi, \tau)^{-1}  r\# \delta_{\tau}, 
\quad r\in R^{\chi}; \qquad \psi (R_{g}^{\chi}) = R^{\psi}_{\Psi(\chi, g)}.$$ 
 See \cite[Lemma 2.3]{AS2}.
Note that \eqref{2.7} is now $\psi(rs) = \omega(\chi, \tau)\psi(r)\psi(s)$,
$r\in R^{\chi}$, $s\in R^{\tau}$.
\begin{lema}\label{c-adj} If $r\in R^{\chi}_{g}$ and $s\in R^{\tau}$ then
\begin{equation} 
 \psi\left( [r, s]_{c}\right) =  
\omega(\chi, \tau) [\psi(r), \psi(s)]_{c}.\end{equation}
\end{lema}
\pf We have 
\begin{align*}  \psi\left( [r, s]_{c}\right) &= 
\psi\left( rs - \tau(g) sr\right) 
\\ &=  \omega(\chi, \tau) \psi(r) \psi(s) - \omega(\tau, \chi)\tau(g)
\psi(s) \psi(r) \\
 &=  \omega(\chi, \tau) \left(\psi(r) \psi(s) - \langle \tau,\Psi(\chi, g)\rangle
\psi(s) \psi(r) \right) \\
&=  \omega(\chi, \tau) [\psi(r)),  \psi(s)]_{c},
\end{align*}where we used \eqref{2.81}. \epf

\begin{obs} \label{c-psi}  It is possible to show that
$(\psi \otimes \psi) c(r \otimes s) = F. c_{F}(\psi(r) \otimes\psi(s))$,
for all $r\in R^{\chi}_{g}$, $s\in R^{\tau}$.
\end{obs}

From the previous considerations and Lemma \ref{propspsi} we immediately get

\begin{prop}\label{twisted}Let $R$ be an algebra in $_{\Gamma}^{\Gamma}\mathcal{YD}$,
$(x_{i})_{i\in I}$ a family of elements of $R$, $x_{i} \in R_{g_{i}}^{\chi_{i}}$
for some $g_{i}\in \Gamma$, $\chi_{i}\in \VGamma$. Then:

(i). $\psi(\ku\langle x_{i}\rangle) = \ku\langle \psi(x_{i})\rangle$, $\psi(\langle x_{i}\rangle) = \langle \psi(x_{i})\rangle$.

(ii). If $R$ has a presentation by generators $x_{i}$ and relations $t_{j}$,
where also the $t_{j}$'s are homogeneous then $R_{F}$ 
has a presentation by generators $\psi(x_{i})$ and relations $\psi(t_{j})$.

(iii). If $x_{i}$ is central and $\omega(\chi_i, \tau) = \omega(\tau, \chi_i)$ for 
all $\tau$ such that $R^{\tau}\neq 0$,  then  $\psi(x_{i})$ is central.
\qed\end{prop}

\bigbreak
\section{Root vectors and Quantum Serre relations \label{pres-nich}}

\subsection{Root vectors\label{rootvectors}}

In this Section, we fix: 

$\Gamma$ a finite abelian group, $(a_{ij})_{1\le i, j\le \theta}$
a finite Cartan matrix,
$g_{1}, \dots, g_{\theta} \in \Gamma$, $\chi_{1}, \dots, \chi_{\theta}
\in \widehat{\Gamma}$ such that  \eqref{diagnotriv} and \eqref{cartantype} hold. 
Let $d_{1}, \dots, d_{\theta} \in \{1,2,3\}$ such that $d_{i} a_{ij} 
= d_{j}a_{ji}$ for all $i$, $j$.
We set $q_{i} = \chi_{i}(g_{i})$, $N_{i}$  the order of $q_{i}$. 
We assume, for all $i$ and $j$, that the order of $\chi_{i}(g_{j})$ is odd, and 
that $N_i$ is not divisible 
by 3 if $i$ belongs to a connected component of type $G_2$.

\smallbreak
Let $\mathcal X$ be the set of connected components of the Dynkin diagram
corresponding to $(a_{ij})$. 
We assume that for each $I\in \mathcal X$, there exist $c_I, d_I$
such that $I = \{j: c_I\le j \le d_I\}$; that is, after reordering
the Cartan matrix is
 a matrix of blocks corresponding to the connected components.
Let $I\in {\mathcal X}$ and $i \sim j$ in $I$;
then $N_{i} = N_{j}$, hence $N_{I} := N_{i}$ is well defined. 
Let $\Phi_{I}$, resp. $\Phi_{I}^{+}$, be the root system, resp. the subset of positive roots,  
corresponding to the Cartan matrix
$(a_{ij})_{i, j\in I}$; then $\Phi = \bigcup_{I\in \mathcal X}\Phi_{I}$, resp. $\Phi^+ = 
\bigcup_{I\in \mathcal X}\Phi_{I}^{+}$ is the root system, resp. the subset of positive roots,  
corresponding to the Cartan matrix
$(a_{ij})_{1\le i, j\le \theta}.$
 Let $\alpha_{1}, \dots,\alpha_{\theta}$  be the set of simple roots.

\smallbreak
Let $\mathcal W_I$ be the Weyl group  corresponding to the Cartan matrix 
$(a_{ij})_{i, j\in I}$; we identify it with a subgroup of the Weyl group 
$\mathcal W$ corresponding to the Cartan matrix 
$(a_{ij})$. We fix
a reduced decomposition of the longest element $\omega_{0, I}$ of $\mathcal W_{I}$
in terms of simple reflections.
Then we obtain a reduced decomposition of the longest element $\omega_{0}
= s_{i_1} \dots s_{i_P}$ of $\mathcal W$
from the expression of $\omega_{0}$ as product of the $\omega_{0, I}$'s in some
fixed order of the components, say the order arising from the order of the vertices.
Therefore
$\beta_{j} := s_{i_1} \dots s_{i_{j-1}}(\alpha_{i_j})$ 
is a numeration of $\Phi^+$.

\bigbreak
We fix a finite dimensional Yetter-Drinfeld module $V$ over $\Gamma$ 
with a basis $x_1, \dots, x_{\theta}$ with $x_{i}\in V^{\chi_i}_{g_i}$, $1\le i \le \theta$.

\bigbreak
Major examples of modules of Cartan type are the Frobenius-Lusztig kernels.
Let $N > 1$ be an odd natural number and let $q\in \ku$ be a primitive
$N$-th root of 1, not divisible by 3 in case $(a_{ij})$ has a component of type 
$G_2$. Let $\BGamma = \Z/(N)^{\theta} = \langle e_{1}\rangle 
\oplus \dots \oplus \langle e_{\theta}\rangle$; let 
$\eta_{j} \in \bVGamma$ be the unique character
such that $ \langle \eta(j), e(i)\rangle = q^{d_{i}a_{ij}}$. Let
$\mathbb V$ be a Yetter-Drinfeld module over $\BGamma$ with a basis
$X_{1}, \dots, X_{\theta}$ such that
$$ X_{i} \in {\mathbb V}_{\eta_{i}}^{e_{i}}, \qquad 1\le i \le \theta.
$$
We denote by $\mathfrak c$ the braiding of $\mathbb V$.
Lusztig defined root vectors $X_{\alpha} \in {\mathfrak B}({\mathbb V})$,
$\alpha \in \Phi^{+}$ \cite{L2}. One can see from \cite{L3} that,
up to a non-zero scalar, each root vector can be written as an 
iterated braided commutator in some sequence 
$X_{\ell_{1}}, \dots, X_{\ell_{a}}$ of simple root vectors
such as $[[X_{\ell_{1}}, [X_{\ell_{2}}, X_{\ell_{3}}]_{\mathfrak c}]_{\mathfrak c},
[X_{\ell_{4}}, X_{\ell_{5}}]_{\mathfrak c}]_{\mathfrak c}$. 
This can also be seen in the situation in \cite{Ri}.

\bigbreak
We now fix for each $\alpha \in \Phi^{+}$ such a representation
of $X_{\alpha}$ as an iterated braided commutator. In the general case of our $V$,
we define root vectors $x_{\alpha}$ in the tensor algebra $T(V)$,
$\alpha \in \Phi^{+}$, as the same formal iteration of braided commutators
in the elements $x_{1}, \dots, x_{\theta}$ instead of
$X_{1}, \dots, X_{\theta}$ but with respect to the braiding $c$ given by the 
general matrix $\left( \chi_j(g_i)\right)$. Note that each $x_{\alpha}$ is homogeneous
and has the same degree as $X_{\alpha}$, where we mean the degree in the sense
of \cite{L3}. Also,
\begin{equation}\label{bihomog}
x_{\alpha}\in T(V)^{\chi_{\alpha}}_{g_{\alpha}},
\end{equation}
where $g_{\alpha} = g_{1}^{b_1} \dots g_{\theta}^{b_{\theta}}$, $\chi_{\alpha} = 
\chi_{1}^{b_1} \dots \chi_{\theta}^{b_{\theta}}$, where $\alpha = b_1\alpha_1 +
\dots + b_{\theta}\alpha_{\theta}$.

\begin{teo}\label{flk} 
The Nichols algebra ${\mathfrak B}({\mathbb V})$ is presented by generators $X_{i}$, $1\le i \le \theta$, and relations
\begin{align}
\ad_{c}(X_{i})^{1-a_{ij}}(X_{j}) & = 0, \qquad i\neq j, \\
X_{\alpha}^{N}&=0, \qquad \alpha \in \Phi^+.
\end{align}
Moreover, the following elements constitute a basis of 
${\mathfrak B}({\mathbb V})$:
$$X_{\beta_{1}}^{h_{1}} X_{\beta_{2}}^{h_{2}} \dots X_{\beta_{P}}^{h_{P}}, 
\qquad 0 \le h_{j} \le N - 1, \quad  1\le 
j \le P.$$
\end{teo}
\pf It follows from results of Lusztig \cite{L1}, \cite{L2}, 
Rosso \cite{Ro1}, \cite{Ro2} and M\"uller \cite{Mu} that
${\mathfrak B}({\mathbb V})$ is the positive part of the so-called
Frobenius-Lusztig kernel corresponding to the Cartan matrix $(a_{ij})$. 
See \cite[Th. 3.1]{AS2} for details. 
The presentation by generators and relations 
follows from the considerations in the last paragraph of p. 15 and
the first paragraph of p. 16 in \cite{AJS} referring to 
\cite[\S 19, Corollary in p. 120]{dCP}.
The statement about the basis is \cite{L1, L2}.
\epf

\subsection{Nichols algebras of Cartan type\label{dos.4}}

We can now prove the first main result of the present paper, describing 
${\mathfrak B}(V)$  by generators  and relations when $V$ is of finite Cartan type,
improving \cite[Th. 1.1 (i)]{AS2}. As in {\it loc. cit.}, we use repeatedly Remark
\ref{nodep}.

\begin{teo}\label{flk-tw}  
The Nichols algebra ${\mathfrak B}(V)$ is presented by generators 
$x_{i}$, $1\le i \le \theta$, and relations
\begin{align} \label{serrebis}
\ad_{c}(x_{i})^{1-a_{ij}}(x_{j}) & = 0, \qquad i\neq j, \\ \label{rootbis}
x_{\alpha}^{N_{I}}&=0, \qquad \alpha \in \Phi_{I}^+, I\in \mathcal X.
\end{align}
Moreover, the following elements constitute a basis of 
${\mathfrak B}(V)$:
$$x_{\beta_{1}}^{h_{1}} x_{\beta_{2}}^{h_{2}} \dots x_{\beta_{P}}^{h_{P}}, \qquad 0 \le h_{j} \le N_{I} - 1, \text{ if }\, \beta_j \in I, \quad  1\le 
j \le P.$$
\end{teo}

\pf (a) Let us first assume that the braiding is symmetric, that is $\chi_{i}(g_{j}) 
= \chi_{j}(g_{i})$ for all $i,j$.
By \cite[Lemma 4.2]{AS2} we can assume moreover that the Cartan matrix $(a_{ij})$ is connected. 
From our assumptions on the orders of the $\chi_{i}(g_{j})$ we then conclude that the braiding has 
the form $\chi_{j}(g_{i}) = q^{d_{i}a_{ij}}$  for all $i,j$
 where $q$ is a root of unity of order $N = \chi_{i}(g_{i})$. 
See \cite[Lemma 4.3]{AS2}.
Hence the Theorem follows directly from Theorem \ref{flk} and  Remark \ref{nodep}.

(b) In the case of an arbitrary braiding we know from  Lemma 4.1 
of \cite{AS2} that there exists a finite abelian group $\BGamma$ satisfying: 
\begin{itemize}
\item The braiding $c$ 
of $V$ can be realized from a Yetter-Drinfeld module structure over 
$\BGamma$ that we continue denoting by $V$, {\it cf.} Remark
\ref{nodep}.
\item There exists a cocycle $\omega: \bVGamma \times \bVGamma \to 
\ku^{\times}$
with corresponding $F \in \ku\BGamma \otimes \ku \BGamma$ 
such that the braiding of $V_{F}$ is symmetric. Let $\psi: {\mathfrak B}(V) \to {\mathfrak B}(V_F)$ 
be the isomorphism having the same meaning as in \eqref{defpsi}. 
\item The braiding of $V_{F}$ is given in the  basis $\psi(x_i) \otimes \psi(x_j)$ by 
a matrix $(b_{ij}^F)$ such that $b_{ii}^F = \chi_{i}(g_{i})$ and the order of $(b_{ij}^F)$
is again odd for all $i$ and $j$.
\end{itemize}
If $\varpi: T(V) \to {\mathfrak B}(V)$, $\varpi_F: T(V_F) \to {\mathfrak B}(V_F)$ denote 
the canonical maps, then we have a commutative diagram
$$\begin{CD}
T(V) @>{\varpi}>> {\mathfrak B}(V) \\
@V{\psi}VV
@VV{\psi}V \\
T(V_F) @>{\varpi_F}>> {\mathfrak B}(V_F).
\end{CD}$$
Clearly, $\psi(\ker \varpi) = \ker \varpi_{F}$; if $(r_j)_{j \in J}$
is a set of generators of the ideal $\ker \varpi$ with
$r_j \in T(V)^{\eta_j}_{h_j}$ then by Proposition \ref{twisted} $(\psi(r_j))_{j \in J}$
is a set of generators of the ideal $\ker \varpi_{F}$. 
By the symmetric case (a), we know the generators of $\ker \varphi_{F}$.
Let us denote 
$X_{i} := \psi(x_{i})$. Then by Lemma \ref{c-adj} and \eqref{2.7}, we have
$\psi\left(\ad_{c}(x_{i})^{1-a_{ij}}(x_{j})
\right) = u_{ij}\ad_{c}(X_{i})^{1-a_{ij}}(X_{j})$ and
$\psi\left(x_{\alpha}^{N_I}\right) = u_{\alpha}X_{\alpha}^{N_I}$, $\alpha \in \Phi_I^+$ where $u_{ij}, u_{\alpha}$ are non-zero scalars. 
This implies the first claim of the Theorem. 
The second follows in a similar way.   
\epf

Let $\wtoba(V)$ be the braided Hopf algebra in $\ydg$ generated by 
$x_1, \dots, x_{\theta}$ with relations \eqref{serrebis}, where the $x_{i}$'s are 
primitive. Let $\ftoba(V)$ be the subalgebra of $\wtoba(V)$ generated by 
$x_{\alpha}^{N_{I}}$, $\alpha\in \Phi_I^+$, $I\in \mathcal X$; it is a Yetter-Drinfeld 
submodule of $\wtoba(V)$.

\begin{teo}\label{powrootvec-subalg} $\ftoba(V)$  is a 
braided Hopf subalgebra in $\ydg$ of $\wtoba(V)$. \end{teo}

\pf (a). As in the proof of Theorem \ref{flk-tw} we first assume that the braiding 
is symmetric. If $i\neq j$, then $\chi_j(g_i)\chi_i(g_j) = 1$ and hence the corresponding 
Serre relation  \eqref{serrebis} says that $x_{i}x_{j} = x_{j}x_{i}$. Thus, we can easily
reduce to the connected case. In such case, $\chi_{j}(g_{i}) = q^{d_{i}a_{ij}}$  as before and the Theorem is shown in \cite{dCP}.

(b). In the general case, we change the group as in the proof of Theorem \ref{flk-tw}. 
The isomorphism $\psi: T(V) \to T(V_F)$ respects the Serre relations up to non-zero scalars
by Lemma \ref{c-adj}. Also, it maps subcoalgebras stable under the action of the group
to subcoalgebras by Lemma \ref{propspsi} (iii). We conclude from (a) that $\ftoba(V)$
is a subcoalgebra of $\wtoba(V)$.  \epf

\section{Linking datum and glueing of connected components}\label{link-dat}
\subsection{Linking datum}\label{linkdat}

In this Section, we fix: $\Gamma$ a finite abelian group, $(a_{ij})_{1\le i, j\le \theta}$
a finite Cartan matrix,
$g_{1}, \dots, g_{\theta} \in \Gamma$, $\chi_{1}, \dots, \chi_{\theta}
\in \widehat{\Gamma}$ such that  \eqref{diagnotriv} and \eqref{cartantype} hold. 
We preserve the conventions and hypotheses from Section \ref{pres-nich}.

\begin{defi}We say that two vertices $i$ and $j$  
{\em are linkable} (or that $i$ {\em is linkable to} $j$) if

\begin{flalign}\label{link0} &i\not\sim j,&\\
\label{link1} &g_{i}g_{j} \neq 1 \text {  and}& \\
\label{link2}
&\chi_{i}\chi_{j} = 1.&
\end{flalign}

If $i$ is linkable to $j$, then  $\chi_{i}(g_{j})\chi_{j}(g_{i}) = 1$ by \eqref{link0};
it follows then from \eqref{link2} that
\begin{equation}\label{link3} \chi_{j}(g_{j}) = \chi_{i}(g_{i})^{-1}.
\end{equation}

\begin{lema}\label{sndproplink}Assume that $i$ and $k$, resp. $j$ and $\ell$, are linkable. 
Then $a_{ij} = a_{k\ell}$, $a_{ji} = a_{\ell k}$. In particular, a vertex $i$ can not be linkable to two different vertices $j$ and $h$.
\end{lema}
\pf If $a_{i\ell} \neq 0$ then $a_{ij} = a_{ji} = 0$ (otherwise $j\sim \ell$) and
$a_{k\ell} = a_{\ell k} = 0$ (otherwise $i\sim k$). If $a_{jk} \neq 0$ then $a_{ij} = a_{ji} = 0$  (otherwise $i\sim k$) and
$a_{k\ell} = a_{\ell k} = 0$ (otherwise $j\sim \ell$).
Assume that $a_{i\ell} = 0 = a_{jk}$. Then
$$
\chi_{i}(g_{i})^{a_{ij}} = \chi_{i}(g_{j})\chi_{j}(g_{i}) =
\chi_{k}^{-1}(g_{j}) \chi_{\ell}^{-1}(g_{i}) =
\chi_{j}(g_{k})\chi_{i}(g_{\ell}) =
\chi_{\ell}^{-1}(g_{k}) \chi_{k}^{-1}(g_{\ell}) =
\chi_{k}(g_{k})^{-a_{k\ell}} =\chi_{i}(g_{i})^{a_{k\ell}}.
$$
Then $N_i$ divides $a_{ij} - a_{k\ell}$ and analogously, $N_k$ divides $a_{ij} - a_{k\ell}$.
So that $a_{ij} = a_{k\ell}$ by the assumptions on the order of $N_i$ and $N_k$; by symmetry, $a_{ji} = a_{\ell k}$. Assume that
a vertex $i$ is linkable to  $j$ and $h$. Then $2 = a_{ii} = a_{jh}$, so $j=h$.
\epf

A {\em linking datum}  for $\Gamma$, $(a_{ij})$,
$g_{1}, \dots, g_{\theta}$ and $\chi_{1}, \dots, \chi_{\theta}$ is 
a collection $(\lambda_{ij})_{1 \le i < j \le \theta, \, i\not\sim j}$ of elements in $\ku$ 
such that $\lambda_{ij}$ is arbitrary if $i$ and $j$  
are linkable but 0 otherwise. Given a linking datum,
we say that two vertices $i$ and $j$  
{\em are linked} if $\lambda_{ij}\neq 0$.
\end{defi}

This definition generalizes part of the definition of compatible datum in \cite[Section 5]{AS1}. 
We shall represent a linking datum by the Dynkin diagram of the Cartan matrix $(a_{ij})$
joining linked vertices by a dotted line. To have a complete picture we add 
the pair $(g_i, \chi_i)$ below the vertex $i$.

\begin{defi}\label{abodrio} Let us fix a decomposition $\Gamma = <Y_1>\oplus \dots \oplus <Y_s>$; let $M_h$ denote the order of $Y_h$, $1\le h \le s$. 
We denote by $\abodrio$ the algebra presented by generators $a_{1}, \dots, a_{\theta}$, 
$y_{1}, \dots, y_{s}$ and relations 
\begin{equation}
\label{relations43}
y_{h}^{M_h} = 1,  
\quad y_{m}y_{h} = y_{h}y_{m}, \quad 1 \le m,h  \le s,\end{equation}
\eqref{relations1}, 
\eqref{relations2},
\eqref{relations3} and 

\begin{equation}
\label{relations44}
a_{\alpha}^{N_{I}} = 0,  
\qquad \alpha \in \Phi_I^{+}, I\in \mathcal X.
\end{equation}
 \end{defi}

\begin{obs}\label{link0o1}
In the preceding definition, one could consider only linking data with $\lambda_{ij} = 
1$ or 0. Indeed, one can replace the generator $a_{i}$ by $\lambda_{ij}^{-1}a_{i}$ whenever
$\lambda_{ij} \neq 0$ for some $j$ which is unique by Lemma \ref{sndproplink}. The other relations do not change since they are homogeneous in the $a_{i}$'s. However, in
the more general case where the relations \eqref{relations44} have a non-zero right side, one needs general linking data.
\end{obs}

\begin{exa}Here is a linking datum where all the connected components are points:
$$
\matrix \bullet & \hdots  & \bullet\\
\bullet & \quad  & \bullet\\
\bullet & \hdots  & \bullet \endmatrix
$$
\end{exa}

\begin{exa}\label{FLK} Let ${\bf B} :=(b_{ij})_{1\le i, j\le R}$ 
be a finite Cartan matrix, $0\le M \le R$ and 
$q\in \ku$ a root of unity of  order $N$; we assume $N$ is odd, and prime to 3 if 
${\bf B}$ contains a component of type $G_2$. Let
$d_1, \dots, d_R$ be integers in $\{1, 2, 3\}$ such that $d_ib_{ij} =
d_jb_{ji}$.
Let $\theta = R + M$, $\widetilde{\bf B} :=(b_{ij})_{1\le i, j\le M}$ and 
${\bf A} = (a_{ij})$ be the Cartan matrix
$$
{\bf A} = \pmatrix {\bf B} &  0\\
0& \widetilde{\bf B}\endpmatrix. $$
Let $\Gamma = (\Z/(N))^R$, $g_{1}, \dots, g_{R}$ the canonical basis of $\Gamma$
and 
$\chi_{1}, \dots, \chi_{R}$ be the character given by $\chi_i(g_j) = q^{d_ib_{ij}}$;
let  $g_{R+ j} = g_{j}$, $\chi_{R+ j} = \chi_{j}^{-1}$, $1\le  j\le M$.
Note that $j$ and $j+R$ are linkable, $1\le  j\le M$.
Finally, let $\lambda_{j, j+R} = 1$ if $1\le  j\le M$ and 0 otherwise;
then $(\lambda_{ij})_{1 \le i < j \le \theta}$ is a linking datum  
for $\Gamma$, $(a_{ij})$,
$g_{1}, \dots, g_{\theta}$ and $\chi_{1}, \dots, \chi_{\theta}$.
The Hopf algebra $\abodrio$
with comultiplication determined  by \eqref{sk-gl} 
is the parabolic part of a Frobenius-Lusztig kernel.
Since the numeration of the Dynkin diagram is so far arbitrary,
any such parabolic appears in this way.
\end{exa}

\begin{exa}Here are some exotic examples of linking data:

 Take 4 copies of $A_3$ and label the vertices such that $\{1,2,3\}$,
$\{4,5,6\}$, $\{7,8,9\}$ and $\{10,11,12\}$ are the connected components. Then link
3 with 4, 6 with 7, 9 with 10 and 12 with 1. It is possible to realize this linking
over $\Z/(N)^{12}$ for any odd $N$; the corresponding braiding will  be 
symmetric in each component, that is, the corresponding subalgebra is the "Borel" part
of a Frobenius-Lusztig kernel. More examples arise considering more copies of
more general components.
\end{exa}

\subsection{Altering the multiplication by a cocycle}\label{altering}
The following variation of Drinfeld's twisting was stated by Doi:
if $H$ is a Hopf algebra and $\sigma: H\times H \to \ku$ is an invertible 2-cocycle,
so that
\begin{align*} \sigma(x_{(1)}, y_{(1)}) \sigma(x_{(2)} y_{(2)}, z) &=
\sigma(y_{(1)}, z_{(1)}) \sigma(x, y_{(2)}z_{(2)}), \\
\sigma(1, 1) &=1,
\end{align*}
for all $x,y,z\in H$, then $H_{\sigma}$-- the same $H$ but with the multiplication
$._{\sigma}$ below-- is again a Hopf algebra, where
$$
x._{\sigma}y = \sigma(x_{(1)}, y_{(1)}) x_{(2)} y_{(2)} \sigma^{-1}(x_{(3)}, y_{(3)}).
$$
\begin{lema}\label{doitak}\cite{DT} Let $U$, $B$ be Hopf algebras.

(a) Let $\tau:U \otimes B \to \ku$ be a  bilinear map such that
for all $u, v\in U$, $a,b\in B$
\begin{enumerate}
\item $\tau(uv, a) = \tau(u, a_{(1)})\tau(v, a_{(2)})$,
\item $\tau(u, ab) = \tau(u_{(1)}, b)\tau(u_{(2)}, a)$,
\item $\tau(1, a) = \varepsilon(a)$,
\item $\tau(u, 1) = \varepsilon(u)$.
\end{enumerate}
Let $H$ be the tensor product Hopf algebra $H = U \otimes B$ and let 
$\sigma: H\otimes H \to \ku$ be the bilinear map $\sigma(u\otimes a, v\otimes b)
= \varepsilon(u) \tau(v, a)\varepsilon(b)$, for all
$u, v\in U$, $a,b\in B$. Then $\tau$ is convolution invertible with inverse given by $\tau^{-1}(v, a)= \varphi(\Ss v)(a) =  \varphi(v)(\Ss^{-1} a)$;
$\sigma$
is an invertible 2-cocycle-- with inverse $\sigma^{-1}(u\otimes a, v\otimes b)
= \varepsilon(u) \tau^{-1}(v, a)\varepsilon(b)$, for all
$u, v\in U$, $a,b\in B$-- and consequently $H_{\sigma}$ is a Hopf algebra.

\bigbreak
(b)  Assume that $B$ is finite dimensional and let $\varphi: U \to (B^*)^{{\rm cop}}$
be a Hopf algebra homomorphism. Then $\tau:U \otimes B \to \ku$, $\tau(v, a)= \varphi(v)(a)$,
is invertible-- with inverse given by $\tau^{-1}(v, a)= \varphi(\Ss v)(a) =  \varphi(v)(\Ss^{-1} a)$,
and satisfies 1, 2, 3 and 4. Reciprocally, given such $\tau$ there is a unique such $\varphi$.
\qed
\end{lema}

The following result is probably known. We include it for completeness.

\begin{lema}\label{doitak2} Let $U$, $B$ and $\tau$ be as in the preceeding Lemma.
Assume that 
\begin{itemize}
\item $U$ is generated as an algebra by skew-primitive elements $u_i$, $i\in I$ and group-like elements $g_k$, $k\in K$, which in addition generate $G(U)$ as a monoid;
\item $B$ is generated as an algebra by skew-primitive elements $b_j$, $j\in J$ and group-like elements $h_{\ell}$, $\ell\in L$, which in addition generate $G(B)$ as a monoid.
\end{itemize}
Let $A$ be an algebra and let $\alpha: U\to A$, $\beta: B\to A$ be algebra maps and let
$\gamma: (U\otimes B)_{\sigma} \to A$ be given by $\gamma(u\otimes b) = \alpha(u)\beta(b)$
for all $u\in U$, $b\in B$. Then $\gamma$ is an algebra map if and only if 
\begin{equation}\label{salemapa}
\tau(u_{(1)}, b_{(1)}) \alpha(u_{(2)}) \beta(b_{(2)}) = \beta(b_{(1)}) \alpha(u_{(1)})
\tau(u_{(2)}, b_{(2)}),
\end{equation}
whenever $u$, resp. $b$, belongs to the family $u_i$, $i\in I$ or $g_k$, $k\in K$,
resp.  $b_j$, $j\in J$ or $h_{\ell}$, $\ell\in L$.
\end{lema}
\pf (Sketch). Clearly, $\gamma$ is an algebra map if and only if \eqref{salemapa} holds for all
$u\in U$, $b\in B$. It follows also easily that \eqref{salemapa} holds when $u=1$, or $b=1$, or $u\in G(U)$ and $b\in G(B)$. Next, let $u, v\in U$ and $b, c\in B$ be arbitrary elements; one
can then check that \eqref{salemapa} holds for $uv$ and $bc$ if  it holds for
all the possibilities $u_{(1)}$ and $c_{(1)}$; $u_{(2)}$ and $b_{(1)}$; $v_{(1)}$ and $c_{(2)}$;
$v_{(2)}$ and $b_{(2)}$. From this observation and the hypothesis the Lemma follows.
\epf

\bigbreak
\subsection{Glueing of connected components}\label{glueing}
In this subsection, we  fix a linking datum  $(\lambda_{ij})_{1 \le i < j \le \theta, \, i\not\sim j}$ 
for $\Gamma$, $(a_{ij})_{1\le i, j\le \theta}$,
$g_{1}, \dots, g_{\theta}$ and $\chi_{1}, \dots, \chi_{\theta}$. 

\

We denote  ${\mathcal A} := \abodrio$.

\begin{teo}\label{linking} (a) There exists a unique Hopf algebra structure in 
$\mathcal A$ determined by \eqref{sk-gl}.

\

(b) The dimension of $\mathcal A$ is $\vert \Gamma\vert\prod_{I\in \mathcal X}N_{I}^{\vert \Phi_I^+ \vert}$.
\end{teo}
\pf By induction on the number of connected components. Here is the first step:

\begin{lema}\label{linking-conn} 
Theorem \eqref{linking} is true if 
the Dynkin diagram corresponding to $(a_{ij})_{1\le i, j\le \theta}$ is connected.
\end{lema}
\pf Let $V = \oplus_{1\le i\le \theta}V^{\chi_{i}}_{g_i}$ be a Yetter-Drinfeld
module over $\Gamma$ with $\dim V^{\chi_{i}}_{g_i} = 1$ and pick 
$x_i\in V^{\chi_{i}}_{g_i} - 0$. By Theorem \ref{flk-tw} and the formulas for the biproduct,
there exists a unique algebra map ${\mathcal F}: {\mathcal A} \to \toba(V)\# \ku\Gamma$ such that
${\mathcal F}(a_i) = x_i\# 1$, ${\mathcal F}(y_t) = 1\# y_t$. Also, by Theorem \ref{flk-tw}
again, there are algebra maps
${\mathcal G}_1: \toba(V) \to {\mathcal A}$, ${\mathcal G}_2: \ku\Gamma 
\to {\mathcal A}$ such that
${\mathcal G}_1(x_i) = a_i$, ${\mathcal G}_2(y_t) = y_t$. 
Let ${\mathcal G}: \toba(V) \# \ku\Gamma\to {\mathcal A}$,
${\mathcal G}(x\# u) = {\mathcal G}_1(x){\mathcal G}_2(u)$, 
$x\in \toba(V)$, $u\in \ku\Gamma$;
then ${\mathcal G}$ is an algebra map by \eqref{relations1}. 
It is clear now that ${\mathcal F}$ is an isomorphism with 
inverse ${\mathcal G}$; thus $\mathcal A$ is a Hopf algebra and has the desired dimension
by the dimension formula in Theorem \ref{flk-tw}. \epf

For the rest of this proof we assume: there exists $\ttheta < \theta$ such that 
$i\sim j$, resp. $i\not\sim h$, if $1\le i \le \ttheta$ and $1\le j \le \ttheta$, resp.
$\ttheta < h\le \theta$. Let $J = \{1, \dots, \ttheta\} \in \mathcal X$.
Let $\Upsilon := <Z_1>\oplus \dots \oplus <Z_{\ttheta}>$,
where the order of $Z_i$ is the least common multiple of $\ord g_{i}$ and $\ord \chi_{i}$, 
$1\le i \le \ttheta$.  Let $\eta_{j}$ be the unique character of $\Upsilon$ such that
$\eta_{j}(Z_{i}) = \chi_{j}(g_{i})$, $1\le i \le \ttheta$, $1\le j \le \ttheta$.
This is well defined because $\ord g_{i}$ divides $\ord Z_{i}$ for all $i$.
\begin{itemize}
\item $\bc := {\mathcal A}\left(\Gamma,\, (a_{ij})_{\ttheta <  i, j\le \theta},\, 
(g_{i})_{\ttheta < i \le \theta}, \, (\chi_{j})_{\ttheta < j \le \theta},  
\, (\lambda_{ij})_{\ttheta < i < j \le \theta,  i\not\sim j} \right)$, 
with generators $b_{\ttheta +1}$, \dots, $b_{\theta}$ (instead of the $a_{i}$'s) and  $y_{1}, \dots, y_{s}$;
\item $\uc := {\mathcal A}\left(\Upsilon,\, (a_{ij})_{1\le i, j\le \ttheta},\, 
(Z_{i})_{1 \le i \le \ttheta}, \, (\eta_{j})_{1 \le j \le \ttheta},  
\, (\lambda_{ij})_{1 \le i < j \le \ttheta, \, i\not\sim j} \right)$, with generators $u_{1}, 
\dots, u_{\ttheta}$ (instead of the $a_{i}$'s) and  $z_1, \dots, z_{\ttheta}$.
\end{itemize}

Note that the linking datum of $\uc$ is empty since $J$ is connected.
By the recurrence hypothesis,  
$\dim \bc = \vert \Gamma\vert \prod_{I\in \mathcal X, I\neq J}N_{I}^{\vert \Phi_I^+ \vert}$
and $\dim \uc = \vert \Upsilon\vert N_{J}^{\vert \Phi_J^+ \vert}$.

\begin{lema}\label{exist-phi} (a). For each $i$, $1 \le i \le \ttheta$,
there exists a unique character $\gamma_{i}: \bc \to \ku$ such that
\begin{equation}\label{gammai}
\gamma_{i}(y_k) = \chi_{i}(y_k), \qquad \gamma_{i}(b_j) = 0,
\end{equation}
$1 \le k \le s$, $\ttheta +1 \le j \le \theta$.

(b). Let $({\widetilde \lambda}_{ij})_{1 \le i < j \le \theta, \, i\not\sim j}$ 
be an arbitrary linking datum. For each $i$, $1 \le i \le \ttheta$,
there exists a unique $(\varepsilon, \gamma_i)$--derivation 
$\delta_{i}: \bc \to \ku$ such that
\begin{equation}
\delta_{i}(y_k) = 0, \qquad \delta_{i}(b_j) = {\widetilde \lambda}_{ij},
\end{equation}
$1 \le k \le s$, $\ttheta +1 \le j \le \theta$.

(c). There exists a unique Hopf algebra map $\varphi: \uc \to (\bc^*)^{{\rm cop}}$ such that
\begin{equation}
\varphi(z_i) = \gamma_{i}, \qquad \varphi(u_i) = \delta_{i},
\end{equation}
$1 \le i \le \ttheta$. \end{lema}

\pf
(a). We have to show that $\gamma_{i}$ preserves the relations \eqref{relations43},
\eqref{relations1}, \eqref{relations2}, \eqref{relations3}, \eqref{relations44}. This is clear for  \eqref{relations43}, \eqref{relations1},
\eqref{relations2}, \eqref{relations44}. We check \eqref{relations3}: 
let $\widetilde{\theta} + 1 \le j, h \le \theta$ 
such that $\lambda_{jh} \neq 0$. Then 
$$
\chi_{i}(g_{j} g_{h}) = \chi_{j}\chi_{h} (g_{i})^{-1} = 1,
$$
by \eqref{link0} and \eqref{link2}. So that relations \eqref{relations3} hold and (a) is proven.  

(b). This is equivalent to: there exists an algebra map $T: \bc\to M_{2}(k)$ such that
\begin{equation}\label{derivacion}
T(y_{k}) = \pmatrix 1 &0 \\ 0 & 
\gamma_{i}(y_{k})\endpmatrix,  \qquad T(b_{j}) = \pmatrix 0 &
{\widetilde \lambda}_{ij} \\ 0 &  0\endpmatrix,
\end{equation}
$1\le k \le s$, $\widetilde{\theta} + 1 \le j \le \theta$. Then $T$ is  of the
form $T(a) = \pmatrix \epsilon(a) & \delta_{i}(a) \\ 0 & 
\gamma_{i}(a)\endpmatrix$ and $\delta_{i}$ is the desired derivation. So, we need to show 
that the relations  \eqref{relations43},
\eqref{relations1}, \eqref{relations2}, \eqref{relations3}, \eqref{relations44} hold for 
the matrices in \eqref{derivacion}. This is evident for 
\eqref{relations43}.
For \eqref{relations1} it amounts to ${\widetilde \lambda}_{ij} \chi_{i}(y_{k})^{-1} = {\widetilde \lambda}_{ij} \chi_{j}(y_{k})$, which follows from \eqref{link2} when 
${\widetilde \lambda}_{ij} \neq 0$. For \eqref{relations2} and \eqref{relations44} the argument is clear. 
Finally, the left hand side of \eqref{relations3} for $j < h$, is 0,
whereas the right-hand side also vanishes since ${ \lambda}_{jh} \gamma_{i}(g_{j})\gamma_{i}
(g_{h}) = { \lambda}_{jh} \chi_{j}(g_{i})^{-1} \chi_{h}(g_{i})^{-1} = { \lambda}_{jh}$ 
by \eqref{link2} again.

(c). It is enough to verify that $\delta_{i}$, $\gamma_{i}$ satisfy the defining relations 
\eqref{relations43},
\eqref{relations1}, \eqref{relations2}, \eqref{relations3}, \eqref{relations44} for $\uc$. 
Indeed, this will automatically imply that $\varphi$ is a Hopf algebra map. 
Note that \eqref{relations3} are empty since the Dynkin diagram of $\uc$ is connected.
For
\eqref{relations43}, it is enough to verify that the equalities hold 
when applied to the generators 
$b_{\ttheta +1}, \dots, b_{\theta}$, $y_{1}, \dots, y_{s}$ since both sides are algebra maps. 
This is now not difficult; for instance 
$(\gamma_{m}\gamma_{h}) (b_{j}) = \gamma_{m} (g_{j}) \gamma_{h} (b_{j}) + \gamma_{m} (b_{j}) =
0 = (\gamma_{h}\gamma_{m}) (b_{j}).$ 
The first relations in \eqref{relations43} for $\Upsilon$ hold since 
$\ord \gamma_i$ divides $\ord Z_i$ for all $i$.
For \eqref{relations1} we need again to verify 
only on generators, since both sides are skew-derivations; this verification is in turn 
straightforward.  
The left-hand side of the Serre relations
\eqref{relations2} is a skew-derivation by \cite[Lemma A.1]{AS2}; again we are reduced to 
see that $(\ad \delta_{i})^{1 - a_{ij}}\delta_{j} (b_{h}) = 0 = 
(\ad \delta_{i})^{1 - a_{ij}}\delta_{j}(y_{t})$,
$\ttheta + 1\le h \le \theta$, $1\le t \le s$. Write $(\ad \delta_{i})^{1 - a_{ij}}
\delta_{j} = \delta_{i} (\ad \delta_{i})^{ - a_{ij}}\delta_{j} - \widehat{q} 
\left((\ad \delta_{i})^{- a_{ij}}\delta_{j}\right) \delta_{i}$, where $\widehat{q}$ is a root of 1.
Then 
$(\ad \delta_{i})^{1 - a_{ij}}\delta_{j} (y_{t})  
=\delta_{i} (y_{t}) (\ad \delta_{i})^{ - a_{ij}}\delta_{j}(y_{t}) - \widehat{q} 
(\ad \delta_{i})^{- a_{ij}}\delta_{j} (y_{t}) \delta_{i} (y_{t})=0.
$
Similarly,
$$(\ad \delta_{i})^{1 - a_{ij}}\delta_{j} (b_{h})  
=\delta_{i} (b_{h}) (\ad \delta_{i})^{ - a_{ij}}\delta_{j}(1) - \widehat{q} 
(\ad \delta_{i})^{- a_{ij}}\delta_{j} (g_{h}) \delta_{i} (b_{h})=0,
$$
since $(\ad \delta_{i})^{- a_{ij}}\delta_{j}$ is a homogeneous polynomial in $\delta_{i}$, 
$\delta_{j}$ of positive degree. Finally, relations  \eqref{relations44} follow from the next Lemma. \epf

\begin{lema}\label{potsvan0} Let $B$ be a finite dimensional pointed Hopf algebra  
generated as an algebra by group-like elements and a family $b_j$, $j\in \mathcal J$, of $(h_j, 1)$-primitives, for some $h_j\in G(B)$. 
Let $\widehat{\uc}$ be the algebra presented by generators   $u_{1}, \dots, u_{\ttheta}$  and  $z_1, \dots, z_{\ttheta}$ with exactly the same relations 
as for $\uc$ except for \eqref{relations44}; it is a Hopf algebra via \eqref{sk-gl}.
Let $N = N_J$.  Assume there exists a Hopf algebra map $\phi: \widehat{\uc} \to (B^*)^{{\rm cop}}$ such that $\gamma_i := \phi(z_{i})$
and $\delta_i := \phi(u_{i})$ satisfy 
\begin{equation}\label{condiciones}
\gamma_i(b_{j})=0,  \qquad \delta_i(g)=0, \qquad j\in \mathcal J,\quad  g \in G(B),
\end{equation}
for all $1\le i \le \ttheta$. Then $\phi(u_{\alpha}^{N}) = 0$ for all $\alpha \in \Phi_J^{+}$.
\end{lema}
\pf There exists a Hopf algebra projection $\varpi: \widehat{\uc} \to \ku \Upsilon$
such that $\varpi(u_i) = 0$ and $\varpi(z_i) = z_i$ for all $i$.
Let $\kc$ be the subalgebra of  $\widehat{\uc}$ generated by $u_{\alpha}^{N}$,
$\alpha \in \Phi_J^{+}$, and $z_i^N$, $1\le i \le \ttheta$.
We claim that 
\begin{equation}\label{proy}
\phi(u) = \phi(\varpi(u))\end{equation} for all $u\in \kc$. Clearly, this implies the Lemma.
By Theorem \ref{powrootvec-subalg}, we know that $\kc$ is a Hopf subalgebra of $\widehat{\uc}$.
We have to prove that $\phi(u) (b) = \phi(\varpi(u)) (b)$ for $b$ a monomial in
the group-likes of $B$ and the $b_j$'s. We do this by induction on the length of the monomial.

We first check the case of length 1.
Here we show more generally 
that $$\phi(u)(g) = \phi(\varpi(u))(g) \quad \text{and} \quad \phi(u)(b_j) =
\phi(\varpi(u))(b_j)$$ for all $g\in G(B)$, $j\in \mathcal J$ and $u\in G(\widehat{\uc})$
or  $u$  of the form $u_{i_1}\dots u_{i_t}z$,
with $z$ group-like, $1\le i_1, \dots, i_t \le \ttheta$, and $t\ge 2$.
Note that each element in $\kc$ is a linear combination of such $u$'s since $N\ge 2$. 
The case when $u$ is a group-like is clear. Let $u = u_{i_1}\dots u_{i_t}z$,
with $z$ group-like, $1\le i_1, \dots,   i_t \le \ttheta$, and $t\ge 2$.
Then $\phi(u) = \delta_{i_1}\dots \delta_{i_t}\phi(z)$, and
\begin{align*}
\phi(u)(g)  &=\delta_{i_1}(g)\dots \delta_{i_t}(g)\phi(z)(g) = 0, \\
\phi(u)(b_j) &= \sum_{1\le r \le t+1} \delta_{i_1}(g_j)\dots \delta_{i_{r-1}}(g_j)
\delta_{i_{r}}(b_j)\delta_{i_{r+1}}(1)\dots \phi(z)(1) = 0,\end{align*}
where we used \eqref{condiciones} and $t\ge 2$. 

Assume then that $b= cd$ where $c$ and $d$
are monomials satisfying the claim.
Since $\phi$ and $\varpi$ are Hopf algebra maps, we have $\phi(u)(cd) =
\phi(u_{(2)})(c) \phi(u_{(1)})(d) =
\phi(\varpi(u_{(2)}))(c) \phi(\varpi(u_{(1)}))(d) = \phi(\varpi(u))(cd)$.  
\epf

We are ready now to conclude the proof of the Theorem. Consider the cocycle $\sigma: 
(\uc \otimes \bc) \otimes (\uc \otimes \bc) \to \ku$ obtained 
 as in Lemma \ref{doitak} from the map $\varphi$ constructed
in Lemma \ref{exist-phi}. Consider the Hopf algebra 
$(\uc \otimes \bc)_{\sigma}$; it has dimension $\vert \Upsilon\vert 
\vert \Gamma\vert\prod_{I\in \mathcal X}N_{I}^{\vert \Phi_I^+ \vert}$.
We claim that the group-like elements $z_{i} \otimes g_{i}^{-1}$ are central in 
$(\uc \otimes \bc)_{\sigma}$ for all $i$. By definition of $(\uc \otimes \bc)_{\sigma}$,
we have to show  for all $u\in \uc$, $b\in \bc$ and $1\le i\le\ttheta$
\begin{equation}\label{centralgrlk}
\varphi(z_{i}) (b_{(1)}) uz_{i} \otimes b_{(2)} g_{i}^{-1} \varphi(z_{i}^{-1}) (b_{(3)})
= \varphi(u_{(1)}) (g_{i}^{-1}) z_{i}u_{(2)} \otimes g_{i}^{-1} b \varphi(\Ss u_{(3)}) 
(g_{i}^{-1}).
\end{equation}
Since $u \otimes b = (u \otimes 1)(1 \otimes b)$ in  $(\uc \otimes \bc)_{\sigma}$ for all $u\in \uc$, $b\in \bc$, it is enough to check \eqref{centralgrlk}
on generators of $\uc$ and $\bc$. This in turn follows easily from the definitions.

\smallbreak
Let $\acw$ be the quotient of $(\uc \otimes \bc)_{\sigma}$ by the central Hopf subalgebra
$\ku[z_{i} \otimes g_{i}^{-1}: 1\le i \le \ttheta]$ with quotient map $\pi$. 
Then $\dim \acw = \vert \Gamma\vert\prod_{I\in \mathcal X}N_{I}^{\vert \Phi_I^+ \vert}$
 by a result of the second author \cite[Th. 3.3.1]{Mo}. Next we claim the existence of a surjective 
algebra map $\mathcal F: \mathcal A \to \acw$ such that 
$${\mathcal F}(a_i) = \pi(u_i\otimes 1), \quad
{\mathcal F}(a_j) = \pi(1\otimes b_j), \quad {\mathcal F}(y_k) = \pi(1 \otimes y_{k}),$$
for $1\le i \le \ttheta$, $\ttheta + 1\le j \le \theta$, $1\le k \le s$.
Again we have to verify the relations \eqref{relations43},
\eqref{relations1}, \eqref{relations2}, \eqref{relations3}, \eqref{relations44}.
Up to \eqref{relations3} these relations already hold in $(\uc \otimes \bc)_{\sigma}$.
For \eqref{relations3}, it is enough to show that 
$$
\pi(u_i\otimes 1)\pi(1\otimes b_j) - \chi_{j}(g_{i}) \pi(1\otimes b_j)\pi(u_i\otimes 1)
= \lambda_{ij} \left(1 - \pi(1\otimes g_{i}g{_j})\right), \qquad 1\le i \le \ttheta, 
\quad \ttheta + 1\le j \le \theta.
$$
A tedious computation shows that the left-hand side is equal to 
$\chi_{j}(g_{i}){\widetilde\lambda}_{ij}\left(\pi(z_{i}\otimes g{_j}) - 1\right)$.
Since $\pi(z_{i}\otimes g{_i}^{-1}) = 1$, we have $\pi(z_{i}\otimes g{_j}) = 
\pi(1\otimes g_{i}g{_j})$. Hence the claim follows if we choose 
${\widetilde\lambda}_{ij} = - \chi_{i}(g_{j}) \lambda_{ij}$
for all $1\le i \le \ttheta$, $\ttheta + 1\le j \le \theta$.

\smallbreak
On the other hand, we have algebra maps $\mathcal G_{1}: \uc  \to \mathcal A$,
$\mathcal G_{2}: \bc  \to \mathcal A$ given by 
${\mathcal G_{1}}(u_i) = a_i$, 
${\mathcal G_{1}}(z_{i}) = g_i$, ${\mathcal G_{2}}(b_j) = a_j$
${\mathcal G_{2}}(y_k) = y_k$, $1\le i \le \ttheta$, $\ttheta + 1\le j \le \theta$, 
$1\le k \le s$. 
Here we use that $\ord g_{i}$ divides $\ord Z_{i}$ for all $i$.
Let $\mathcal  G: \uc \otimes \bc  \to \mathcal A$ be defined by 
$\mathcal  G(u\otimes b) = \mathcal  G_{1} (u)\mathcal  G_{2}(b)$ for all
$u\in \uc$, $b\in \bc$. We claim that $\mathcal  G$ is 
an algebra map. By Lemma \ref{doitak2}, we have to verify
$\varphi(u_{(1)}) (b_{(1)}) \mathcal G_1(u_{(2)}) \mathcal G_2(b_{(2)}) = 
\mathcal G_2(b_{(1)}) \mathcal G_1(u_{(1)})
\varphi(u_{(2)})(b_{(2)})$, for all  generators. 
This is a straightforward task; for the case $u_i$ and $b_j$ we need 
again the condition ${\widetilde\lambda}_{ij} = - \chi_{i}(g_{j}) \lambda_{ij}$.

Since clearly $\mathcal  G$ factorizes 
through $\acw$, $\mathcal F$ is an isomorphism and the Theorem follows. \epf

\begin{obs} Let us consider a Hopf algebra defined as in Definition \ref{abodrio} 
but replacing \eqref{relations44} by \eqref{relations3} {\it only for simple $\alpha$}, and 
demanding \eqref{relations44} for non-simple $\alpha$. Then it is possible to 
prove the analogue of Theorem \ref{linking}. This gives in particular a new proof
of \cite[Proposition 5.2]{AS1} and generalizes \cite[Theorem 3.6]{AS3}. Indeed, the connected components
of the Dynkin diagrams in {\it loc. cit.} are of type $A_1$.
This result was also found independently by A. Masuoka \cite{Ma}.

\end{obs}

\section{Lifting of relations}\label{serre}

In this Section, we assume the situation described in  Section \ref{pres-nich}.
To lift the Serre relations, we need the following Lemma. 

\begin{lema}\label{serre-mild} Let $1\le i \neq j \le \theta$ and let $I$ be the connected
component containing $i$. 

(a). If $i\sim j$, assume that $N_I\neq 3$; if $i\sim j$ and 
$I$ is of type $B_n$, $C_n$ or $F_4$ assume further $N_I\neq 5$.
Then there exists no $\ell$, $1\le  \ell
 \le \theta$, such that $g_{i}^{1-a_{ij}}g_{j} = g_{\ell}$, $\chi_{i}^{1-a_{ij}}\chi_{j} = \chi_{\ell}$.

(b). Assume that $i\sim j$ and $N_I\neq 3$. 
If  $I$ is of type $B_n$, $C_n$ or $F_4$, resp. $G_2$, assume further that $N_I\neq 5$, 
resp. $N_I\neq 7$. Then 
 $\chi_{i}^{1-a_{ij}}\chi_{j} \neq \varepsilon$.

\end{lema}

\pf (a). Assume that $g_{i}^{1-a_{ij}}g_{j} = g_{\ell}$, $\chi_{i}^{1-a_{ij}}\chi_{j} = \chi_{\ell}$ for some $\ell$. Substituting $g_{\ell}$ and $\chi_{\ell}$ in 
$\langle \chi_{\ell}, g_{i}\rangle \langle \chi_{i}, g_{\ell}\rangle = q_i^{a_{i\ell}}$
and using $\langle \chi_{j}, g_{i}\rangle \langle \chi_{i}, g_{j}\rangle = q_i^{a_{ij}}$
we conclude that
\begin{equation}\label{cincouno}
N_i \text{ divides } 2-a_{ij}-a_{i\ell}.
\end{equation}
Changing the r\^oles of $i$ and $j$ we obtain in the same way
\begin{equation}\label{cincodos}
N_j \text{ divides } a_{ji}(1-a_{ij})-a_{j\ell} + 2.
\end{equation}
First assume that $i\not\sim j$. In particular, $a_{ji} = 0 = a_{ij}$
and $a_{i\ell} =0$ or $a_{j\ell} =0$. If $a_{i\ell} =0$, resp. $a_{j\ell} =0$,
then we get from \eqref{cincouno}, resp. \eqref{cincodos}, that $N_i = 2$,
resp. $N_j=2$, which is not possible.

Next assume that $i\sim j$. If $j = \ell$ then $N_i$ divides $2(1-a_{ij})$
by \eqref{cincouno}. The only possibility is $a_{ij} = -2$ and $N_i=3$;
but this was excluded in the hypothesis. If $i = \ell$ then $N_i$ divides $-a_{ij}$
by \eqref{cincouno} and $N_j$ divides $a_{ji}a_{ij} - 2$ by \eqref{cincodos}; 
but this contradicts our general assumptions on the $N_i$'s.

Finally, if $i \neq \ell$ and $j \neq \ell$ then $a_{i\ell} \neq -3$
and $a_{j\ell} \neq -3$. We discuss the different possible values of $a_{ij}$.
If $a_{ij} = 0$ or $-1$, by \eqref{cincouno} and since $N_i$ is odd we see that
$N_i =3$ or 5, cases excluded by hypothesis. If $a_{ij} = -2$ then $a_{ji} = -1$.
By  \eqref{cincodos}, $N_j$ divides $-1-a_{j\ell}$; this discards everything except
$a_{j\ell} = -1$. But in this last case, $a_{i\ell} =0$ and $N_i$ divides 4 by
\eqref{cincouno}, a contradiction. Finally, $a_{ij} = -3$ is impossible by 
analogous arguments.

(b). Assume that $\chi_{i}^{1-a_{ij}}\chi_{j} = \varepsilon$. We consider first the case 
$a_{ij}\neq 0$. Evaluating at $g_i$, we get $q_{i}^{1-a_{ij}} \chi_{j}(g_{i}) = 1$;
hence $q_i = \chi_{i}(g_{j})$. Evaluating at $g_j$, we get then
$q_j = q_{i}^{a_{ij} - 1}$. Since $q_j^{a_{ji}} = q_{i}^{a_{ij}}$ we finally
obtain
\begin{equation}\label{cincotres}
N_i \text{ divides } a_{ij}a_{ji}- a_{ij} - a_{ji}.
\end{equation}
The possible values of $a_{ij}a_{ji}- a_{ij} - a_{ji}$ are 3, 5 or 7,
where 5, resp. 7, is only possible if $I$ is of type $B_n$, $C_n$ or $F_4$, resp. $G_2$.
This contradicts the hypothesis.

We consider finally the case $a_{ij} = 0$; so that $\chi_{i}\chi_{j} = \varepsilon$.
Since $I$ is connected, there is a sequence $i=i_1, i_2, \dots, i_t = j$
of elements in $I$ such that $a_{i_{\ell}i_{\ell +1}}\neq 0, 2$ for all $\ell$, $1\le \ell < t$.
Then
$$
q_{i}^{a_{i_{1}i_{2}} a_{i_{2}i_{3}} \dots a_{i_{t-1}i_{t}}}
= q_{i_2}^{a_{i_{2}i_{1}} a_{i_{2}i_{3}} \dots a_{i_{t-1}i_{t}}}
= \dots = q_{j}^{a_{i_{2}i_{1}} a_{i_{3}i_{2}} \dots a_{i_{t}i_{t-1}}},
$$
by substituting $q_{i}^{a_{i_{1}i_{2}}} = q_{i_2}^{a_{i_{2}i_{1}}}$, then
$q_{i_2}^{a_{i_{2}i_{1}}} = q_{i_3}^{a_{i_{3}i_{2}}}$ and so on.
Note that $q_i = q_j^{-1}$ since $a_{ij} =0$ and $\chi_{i}\chi_{j} = \varepsilon$.
Hence
\begin{equation}\label{cincocuatro}
N_i \text{ divides } a_{i_{1}i_{2}} a_{i_{2}i_{3}} \dots a_{i_{t-1}i_{t}} + 
a_{i_{2}i_{1}} a_{i_{3}i_{2}} \dots a_{i_{t}i_{t-1}}.
\end{equation}
The possible values of the sum in \eqref{cincocuatro} are $\pm 2$ or $\pm 3$.
Hence \eqref{cincocuatro} contradicts our assumptions in (b). \epf

\bigbreak
Let now $A$ be a  pointed Hopf algebra with
$G(A) \simeq  \Gamma$, not necessarily finite dimensional. Let $R$ be the diagram of $A$ (see Section 2.1). We assume there is an 
isomorphism   ${\mathcal P}(R) \cong V$ in $\ydg$. Then
$$ \oplus_{\substack{g,h \in \Gamma \\ \varepsilon \neq \chi \in \VGamma}} 
{\mathcal P}_{g,h}(A)^{\chi} \xrightarrow{\cong} A_{1}/A_{0} \xleftarrow{\cong} V \# \ku \Gamma$$
(see \cite[Lemma 5.4]{AS1}).
Let $a_i \in {\mathcal P}_{g_i,1}(A)^{\chi_i}$, $1\le i \le \theta$,
such that $\overline{a_i}$ is mapped onto $x_i$ for all $i$. 
Then we know from \cite[Lemma 5.4]{AS2} that for all $g\in \Gamma$, $\chi\in \VGamma$
with $\chi\neq \varepsilon$:

\begin{align}\label{lemma54}
{\mathcal P}_{g,1}(A)^{\chi} &\neq 0 \iff \text{ there is some } 1\le \ell\le \theta:
g=g_{\ell}, \chi= \chi_{\ell}; \\ \label{lemma54bis}
{\mathcal P}_{g,1}(A)^{\varepsilon} &= \ku(1-g).
\end{align}

\begin{teo}\label{QSR} Let $A$ and   $a_1, \dots, a_{\theta}$ be as above.

(a). There is a linking datum $(\lambda_{ij})_{1 \le i < j \le \theta, \, i\not\sim j}$
such that \eqref{relations3} holds.

(b). Let $I \in \mathcal X$. Assume that $N_I\neq 3$.
If  $I$ is of type $B_n$, $C_n$ or $F_4$, resp. $G_2$, assume further that $N_I\neq 5$, 
resp. $N_I\neq 7$. Then the quantum Serre relations \eqref{relations2} hold for all
$i\neq j \in I$.
\end{teo}
\pf It is known that $(\ad a_{i})^{1 - a_{ij}}a_{j} \in 
{\mathcal P}_{g_{i}^{1-a_{ij}}g_{j},1}(A)^{\chi_{i}^{1-a_{ij}}\chi_{j}}$, see for instance \cite[Appendix]{AS2}.
Part (b) of the Theorem then follows from Lemma \ref{serre-mild}, \eqref{lemma54} and \eqref{lemma54bis}.

To prove part (a), let us assume that $i \not\sim j$.  By \eqref{lemma54} and \eqref{lemma54bis} again,
 $a_{i}a_{j} - \chi_{j}(g_{i})a_{j}a_{i} = \lambda_{ij}(1 - g_{i}g_{j})$,  
for some $\lambda_{ij}\in \ku$. We can choose 
$\lambda_{ij} = 0$ when $g_{i}g_{j} = 1$ or else if $\chi_{i}\chi_{j} \neq \varepsilon$.
That is,
$\lambda_{ij}$ is a linking datum for $(a_{ij})$, $g_{1}, \dots, g_{\theta}$ 
and $\chi_{1}, \dots, \chi_{\theta}$; and
\eqref{relations3} hold.  \epf

\begin{lema}\label{powrootvec}  Let $A$ and   $a_1, \dots, a_{\theta}$ be as above.
Assume further that 
\begin{itemize}
\item the hypothesis from Theorem \ref{QSR} part (b) holds for all $I\in \mathcal X$.
\item $A$ is   finite dimensional. \item $g_i^{N_i} = 1$, $1\le i\le \theta$. 
\end{itemize}
Then the relations \eqref{relations44} hold in $A$. 
\end{lema}
\pf
Let us fix $I\in \mathcal X$.
Let $\widehat{\uc}$ be the algebra presented by generators ${\widehat a}_{i}$, $i\in I$,
$y_{1}, \dots, y_{s}$ and relations 
\eqref{relations43}, \eqref{relations1}, \eqref{relations2} and \eqref{relations3}; 
it is a Hopf algebra via \eqref{sk-gl}.
Let $N = N_I$ and let $\kc$ be the subalgebra of  $\widehat{\uc}$ generated by 
${\widehat a}_{\alpha}^{N}$,
$\alpha \in \Phi_I^{+}$, and $g_i^N$, $i\in I$. 
By Theorem \ref{powrootvec-subalg}, we know that $\kc$ is a Hopf subalgebra of $\widehat{\uc}$.
Note that $\kc$ is a graded Hopf algebra with trivial coradical.
By the choice of the $a_i$'s in $A$ and Theorem \ref{QSR}, we see there is a well-defined
Hopf algebra map $\widehat{\uc} \to A$ such that ${\widehat a}_{i} \mapsto a_{i}$, $i\in I$.
The image of $\kc$ under this map is a finite dimensional pointed Hopf algebra; it has 
a trivial coradical by \cite{Mo} and therefore it is trivial. This implies the Lemma.
\epf

\begin{teo}\label{genlifting}  Let $A$  be as above and assume that 
\begin{itemize}
\item the hypothesis from Theorem \ref{QSR} part (b) holds for all $I\in \mathcal X$.
\item $\gr A \simeq
\toba(V)  \# \ku \Gamma$, hence $A$ is   finite dimensional. \item $g_i^{N_i} = 1$, $1\le i\le \theta$. 
\end{itemize}
Then there exists a linking datum $(\lambda_{ij})_{1 \le i < j \le \theta, \, i\not\sim j}$
such that $$A \simeq \abodrio.$$
\end{teo}

\pf By Theorem \ref{QSR} and Lemma \ref{powrootvec}, there exists a linking datum
$(\lambda_{ij})_{1 \le i < j \le \theta, \, i\not\sim j}$ and a surjective Hopf algebra map $\mathcal J: \mathcal A \to  A$, where $\mathcal A = \abodrio$.
But 
$\dim A = \dim \gr A = \vert\Gamma\vert\dim \toba(V) = \dim \mathcal A$
by Theorems \ref{flk-tw} and 
\ref{linking}; hence  $\mathcal J$ is an isomorphism.
\epf

\section{Hopf algebras generated in degree one}\label{degree one}

In this Section, $\Gamma$ is a finite abelian group. Let 
$V \in ^{\Gamma}_{\Gamma}\mathcal {YD}$ be of Cartan type with braiding 
$(b_{ij})$ and Cartan matrix $(a_{ij})$. As in \cite{AS2} we say that the braiding 
$(b_{ij})$ satisfies the {\itshape relative primeness condition} if
for all $i,j$,  $(a_{ij})$ is 0 or relatively prime to the order of $b_{ii}$.

The next lemma follows from \cite{AS2}. We will apply it in the case of 
$2 \times 2$ (hence symmetrizable) Cartan matrices.

\begin{lema}\label{finite} 
Let $V \in ^{\Gamma}_{\Gamma}\mathcal {YD}$ be of symmetrizable Cartan 
type with braiding $(b_{ij})_{1 \leq i,j \leq \theta}$. Assume that for 
all $1 \leq i,j \leq \theta$, the order of $b_{ij}$ is odd, and that 
$(b_{ij})_{1 \leq i,j \leq \theta}$ satisfies the relative primeness condition. 
If $\toba(V)$ is finite-dimensional, then $V$ is of finite Cartan type.
\end{lema}
\pf
By \cite[Lemma 4.1]{AS2} we can realize the braiding over a suitable 
finite abelian group
$\widetilde{\Gamma}$ and twist with a 2-cocycle $F$ such that the 
resulting braiding $(b_{ij}^F)$ is symmetric with elements of odd order, 
has the same diagonal elements and is of Cartan type with the same Cartan 
matrix $(a_{ij})$ as $V$. We can now conclude from \cite[Lemma 4.3]{AS2} 
that $(b_{ij}^F)$ is of FL-type (see \cite{AS2}). Let $\widetilde{V}^{F}$ 
be the Yetter-Drinfeld module over $\widetilde{\Gamma}$ with braiding  $(b_{ij}^F)$. 
Since $\toba(V)$ and $\toba(\widetilde{V}^{F})$ have the same dimension, 
$\toba(\widetilde{V}^{F})$ is finite-dimensional. Then $(a_{ij})$ is of 
finite Cartan type by \cite[Theorem 3.1]{AS2}. 
\epf

\begin{lema}\label{Serre2}
Let $S = \oplus_{n \ge 0} S(n)$ be a finite-dimensional graded Hopf algebra in  
$^{\Gamma}_{\Gamma}\mathcal {YD}$ 
such that $S(0) = \ku 1$. Assume that $S(1)$ is of finite Cartan type with basis 
$(x_{i})_{1 \leq i \leq \theta}$, braiding $(b_{ij})_{1 \leq i,j \leq \theta}$
and Cartan matrix $(a_{ij})_{1 \leq i,j \leq \theta}$ as in (\ref{base-cartan})
. For all $1 \leq l \leq \theta$, let $q_{l} = b_{ll}$ and $N_{l} = ord(q_{l})$.

Let $1 \leq i,j \leq \theta$, $i \neq j$, and assume that $N_{i}, N_{j}$ and ord($b_{ij}$) 
are odd, and $N_{i}$ is not divisible by 3 and $ > 7$.

(a) Assume $ i \sim j$ and let $I$ be the connected component containing $i,j$. If the 
type of $I$ is $B_{n}, C_{n} \text{ or } F_{4}$, assume that $N_{i}$ is not divisible 
by 5. If the type is $G_{2}$, assume that $N_{i}$ is not divisible by 5 or 7. 
Then $(\ad_c x_{i})^{1 - a_{ij}} x_{j} = 0$.

(b) Assume $ i \nsim j$ and $q_{i}q_{j} = 1$ or $ord(q_{i}q_{j}) = ord(q_{i})$. 
Then $x_{i}x_{j} - b_{ij} x_{j}x_{i} = 0$.
\end{lema}

\pf
Define $z_1 := x_i$, $z_2 :=  (\ad_c x_{i})^{1 - a_{ij}} x_{j}$. In both cases we 
have to show $z_{2} = 0$. We assume that $z_2$ is not 0. 
Let $g_{i} \in \Gamma, \chi_{i} \in \VGamma$, $1 \leq i,j \leq \theta$, 
with $b_{ij} = \chi_{j}(g_{i})$ for all $i,j$. Then action and coaction on  $z_1, z_{2}$ are given 
by $\delta(z_1) = g_i \otimes z_1$, $\delta(z_2) = g_{i}^{1 - a_{ij}} g_{j}\otimes z_2$ and $h \cdot z_1 
= \chi_{i}(h) z_1$, $h \cdot z_2 = (\chi_{i}^{1 - a _{ij}} \chi_{j})(h) z_2$ for all $h \in \Gamma$. 
Hence $z_1$, $z_2$ are linearly independent since $g_i \neq  g_{i}^{1 - a_{ij}} g_{j}$ 
or $\chi_{i} \neq \chi_{i}^{1 - a _{ij}} \chi_{j}$. (If both equalities would 
hold, then $g_j = g_{i}^{a_{ij}}, \chi_j = \chi_{i}^{a_{ij}}$, and the Cartan 
type condition
$\chi_{i}(g_{j}) \chi_{j}(g_{i}) = \chi_{i}^{a_{ij}}(g_i)$ would give 
$\chi_{i}(g_{i})^{a_{ij}} = 1$, hence $a_{ij} = 0$ and  $g_j = g_{i}^{a_{ij}} = 1$ 
which is impossible.) The braiding $(B_{kl})_{1 \leq k,l \leq 2}$ of the 2-dimensional 
Yetter-Drinfeld module with basis $z_1,z_2$ is given by
\begin{align*}
B_{11} &= \chi_{i}(g_{i}) = q_{i}, &  B_{12} &= (\chi_{i}^{1 - a_{ij}} \chi_{j})(g_{i}) = q_{i}b_{ji}^{-1},\\ 
B_{21} &= \chi_{i}(g_{i}^{1 - a_{ij}}g_{j}) =q_{i}^{1 - a_{ij}}b_{ji},  &      B_{22} &= (\chi_{i}^{1 - a_{ij}}
\chi_{j})(g_{i}^{1 - a_{ij}}g_{j}) = q_{i}^{1 - a_{ij}} q_{j}. 
\end{align*}
Then  $B_{12} B_{21} = q_{i}^{2 - a_{ij}}$. We claim that $(B_{kl})$ is of Cartan type 
and satisfies the relative primeness condition, that is there are integers $A_{12}, A_{21}$ such that 
\begin{align} \label{Cartan12}
q_{i}^{2 - a_{ij}} &= q_{i}^{A_{12}} \text{ ,and } A_{12} \text{ is relatively prime to } N_{i},\\ 
\label{Cartan21}
q_{i}^{2 - a_{ij}} &= (q_{i}^{1 - a_{ij}}q_{j})^{A_{21}}  \text{ ,
and } A_{21} \text{ is relatively prime to ord}(q_{i}^{1 - a_{ij}}q_{j}).
\end{align}
In both cases $2 - a_{ij}$ is relatively prime to $N_{i}$, because of the hypothesis on $N_i$. This shows 
\eqref{Cartan12}.

We now prove \eqref{Cartan21} in case (a). Then $N_{I} = N_{i} = N_{j}$, and it suffices to 
find an integer $A_{21}$ relatively prime to $N_{i}$ with $q_{i}^{2 - a_{ij}} 
= (q_{i}^{1 - a_{ij}}q_{j})^{A_{21}}$. 

First assume that $a_{ij} \neq 0$. Since $a_{ji} \neq 0$ is relatively prime to $N_{i}$, 
it is enough to consider the $a_{ji}$-th power of \eqref{Cartan21}. Since $q_{i}^{a_{ij}} 
= q_{j}^{a_{ji}}$ by the Cartan condition for $(b_{ij})$, we have to solve $(2 - a_{ij})a_{ji} 
\equiv ((1 - a_{ij})a_{ji} + a_{ij})A_{21} \mod{N_{i}}$. Since $(a_{ij})$ is of finite 
Cartan type, the possible values of $(2 - a_{ij})a_{ji}$
are -3, -4, -5, -6, -9 (-4, -6 resp. -5, -9 only occur if the type is $B_{n}, C_{n}$ or 
$F_{4}$ resp. $G_{2}$); the possible values of  $((1 - a_{ij})a_{ji} + a_{ij})$ are -3, 
-5, -7, (-5, resp. -7 only occur if the type is $B_{n}, C_{n}$ or $F_{4}$ resp. $G_{2}$). 
Hence $(2 - a_{ij})a_{ji}$ and $((1 - a_{ij})a_{ji} + a_{ij})$ are relatively prime to $N_{i}$ 
by assumption, and the claim follows. (Note that $a+b$ is never $0$).

If $a_{ij} = 0$, then by connectedness there is a sequence $i=i_1, i_2, \dots, i_t = j$
of elements in $I$ such that $a_{i_{\ell}i_{\ell +1}}\neq 0, 2$
for all $\ell$, $1\le \ell < t$.
Then as in the proof of Lemma \ref{serre-mild} (b),
$$q_{i}^{a} = q_{j}^{b}, \text{ where }a = a_{i_{1}i_{2}} a_{i_{2}i_{3}} \dots 
a_{i_{t-1}i_{t}}, \text{ and } b = a_{i_{2}i_{1}} a_{i_{2}i_{3}} \dots a_{i_{t-1}i_{t}}.$$
Since the possible values of $a,b$ are 1, 2, -1, -2, the $b$-th power of 
\eqref{Cartan21} leads to the congruence $2b \equiv (b + a)A_{21} \mod{N_{i}}$, 
and the claim again follows.

Assume case (b), in particular $a_{ij} = 0$. If $q_{i}q_{j} = 1$, we get a 
contradiction since the algebra generated by $z_{1},z_{2}$ is finite-dimensional, 
hence $B_{22} \neq 1$ by \cite[Lemma 3.1]{AS1}. If $\text{ord}(q_{i}q_{j}) = 
\text{ord}(q_{i})$, \eqref{Cartan21} is solvable since $N_{i}$ is odd.

Thus we have shown that $(B_{kl})$ is of Cartan type and satisfies the relative 
primeness condition. Hence $(B_{kl})$ is of finite Cartan type by Lemma \ref{finite}.
In both cases $A_{12} = 2 - a_{ij} - N_{i}$ is a solution of \eqref{Cartan12}, 
and $- N_{i} < A_{12} \leq 0$. Hence the possible values of $A_{12}$ are 0, -1, 
-2, -3, and we see that $N_{i} \leq 8$. This contradicts our assumption, and we have 
shown the Serre relation $z_{2} = 0$.
\epf

\begin{lema}\label{rootv}
Let $S = \oplus_{n \ge 0} S(n)$ be a finite-dimensional graded Hopf algebra in  
$^{\Gamma}_{\Gamma}\mathcal {YD}$ such that $S(0) = \ku 1$. Assume that $V = S(1)$ 
is of Cartan type with basis $(x_{i})_{1 \leq i,j \leq \theta}$ as described in 
the beginning of Section \ref{pres-nich}. 
Assume the Serre relations $$(\ad_c x_{i})^{1 - a_{ij}} x_{j} = 0 
\text{ for all } 1 \leq i,j \leq \theta, i \neq j \text{ and } i \sim j.$$
Then the root vector relations $$x_{\alpha}^{N_{I}} = 0, \quad\alpha \in \Phi^{+}_{I}, 
\quad I\in \mathcal X,$$ hold in $S$.
\end{lema}

\pf
We fix a connected component $I\in \mathcal X$. Let $V_{I}$ be the Yetter-Drinfeld 
submodule of $V$ with basis $x_i, i \in I$, and $\wtoba(V_{I})$ the quotient of 
$T(V_{I})$ modulo the Serre relations of all elements $x_i, x_j$ with $i \neq j$ in I. 
Let $N_{I} = N$. The map  $\Psi : T(V_{I}) \subset T(V) \to S$ factorizes over 
$\wtoba(V_{I})$, since the Serre relations hold in $S$. By Theorem 
\ref{powrootvec-subalg} the subalgebra $\ftoba(V_{I})$ of $\wtoba(V)$ generated 
by the powers of the root vectors $x_{\alpha}^{N_{I}}$, $\alpha \in \Phi_I^+$, 
is a braided Hopf subalgebra. As a coalgebra, $\ftoba(V_{I})$ is pointed and has 
trivial coradical. Hence $K := \Psi(\ftoba(V_{I}))$ is a finite-dimensional pointed 
and graded  Hopf subalgebra of $S$ in $^{\Gamma}_{\Gamma}\mathcal {YD}$ with trivial 
coradical. We have to show the root vector relation  $x_{\alpha}^{N_{I}} = 0$, 
$\alpha \in \Phi^{+}_{I}$ in $S$, or equivalently that $K$ is one-dimensional, 
that is $ {\mathcal P}(K) = 0$.

Assume ${\mathcal P}(K) \neq 0$. Since ${\mathcal P}(K)$ is in
$^{\Gamma}_{\Gamma}\mathcal {YD}$, there are $ g \in \Gamma, \chi \in \VGamma$ 
with  ${\mathcal P}(K)_{g}^{\chi} \neq 0$. By \cite[Lemma 3.1]{AS1}, we 
conclude $\chi(g) \neq 1$. But this is a contradiction, since for all 
$ g \in \Gamma, \chi \in \VGamma$,  $K_{g}^{\chi} \neq 0$ implies $\chi(g) = 1$. 
For, $K$ is the $\ku$-span of all monomials
 $$\Psi(x_{\beta_{1}}^{N}) \cdots  \Psi(x_{\beta_{m}}^{N}), \: m \geq 1, \: \beta_{1}, \dots, \beta_{m} \in \Phi_{I}^{+}.$$
For any $\beta \in \Phi_{I}^{+}$ there are natural numbers $b_{i}^{\beta}$, $ 1 \leq i \leq \theta$, such that $\beta = \sum_{i = 1}^{\theta} b_{i}^{\beta} \alpha_{i}$, where $\alpha_{1}, \dots, \alpha_{\theta}$ are the simple roots. By \eqref{bihomog}, $x_{\beta} \in T(V_{I})_{g_{\beta}}^{\chi_{\beta}}$, where $g_{\beta} = \prod_{i \in I} g_{i}^{b_{i}^{\beta}}, \chi_{\beta} = \prod_{i \in I} \chi_{i}^{b_{i}^{\beta}}$. Hence for all  $\beta_{1}, \dots, \beta_{m} \in \Phi_{I}^{+}$, 
$$\Psi(x_{\beta_{1}}^{N}) \cdots  \Psi(x_{\beta_{m}}^{N}) \in K_{g}^{\chi}, \quad \mbox{where} \quad  \chi = \chi_{\beta_{1}}^{N} \cdots \chi_{\beta_{m}}^{N},\:  g = g_{\beta_{1}}^{N} \cdots g_{\beta_{m}}^{N}.$$
It remains to show that $\chi(g) = 1$.

Let $\alpha, \beta \in \Phi_{I}^{+}$. Since the braiding is of Cartan type, 
%\begin{equation} \label{character}
%\begin{split}
$$\chi_{\alpha}(g_{\beta}) \chi_{\beta}(g_{\alpha}) = \prod_{i,j} \chi_{i}^{b_{i}^{\alpha}}(g_{j}^{b_{j}^{\beta}}) \prod_{i,j} \chi_{j}^{b_{j}^{\beta}}(g_{i}^{b_{i}^{\alpha}})
 =\prod_{i,j} (\chi_{i}(g_{j}) \chi_{j}(g_{i}))^{b_{i}^{\alpha} b_{j}^{\beta}}
=\prod_{i,j} \chi_{i}(g_{i})^{a_{ij} b_{i}^{\alpha} b_{j}^{\beta}},$$
%\end{split}
%\end{equation}
and
%\begin{equation} \label{character2}
%\begin{split}
%= \prod_{i} %\chi_{i}^{b_{i}^{\alpha}}(g_{i}^{b_{i}^{\alpha}}) \prod _{i < j} %\chi_{i}^{b_{i}^{\alpha}}(g_{j}^{b_{j}^{\alpha}}) %\chi_{j}^{b_{j}^{\alpha}}(g_{i}^{b_{i}^{\alpha}})
$$\chi_{\alpha}(g_{\alpha}) = \prod_{i} \chi_{i}(g_{i})^{b_{i}^{\alpha} b_{i}^{\alpha}} \prod _{i < j} (\chi_{i}(g_{j}) \chi_{j}(g_{i}))^{b_{i}^{\alpha} b_{j}^{\alpha}}
= \prod_{i} \chi_{i}(g_{i})^{b_{i}^{\alpha} b_{i}^{\alpha}} \prod _{i < j} \chi_{i}(g_{i})^{a_{ij} b_{i}^{\alpha} b_{j}^{\alpha}}.$$
%$\end{split}
%$\end{equation}
Hence, since all the $\chi_{i}(g_{i})$ have order N,
$$\chi_{\alpha}^{N}(g_{\beta}^{N}) \chi_{\beta}^{N}(g_{\alpha}^{N}) = 1, \quad 
\chi_{\alpha}^{N}(g_{\alpha}^{N}) = 1.$$
Therefore we obtain
$$\chi(g) = \prod _{i,j} \chi_{\beta_{i}}^{N}(g_{\beta_{j}}^{N}) = \prod _{i}         
\chi_{\beta_{i}}^{N}(g_{\beta_{i}}^{N}) \prod _{i < j} \chi_{\beta_{i}}^{N}(g_{\beta_{j}}^{N})  
\chi_{\beta_{j}}^{N}(g_{\beta_{i}}^{N}) = 1.$$     
\epf

\begin{teo}\label{degree1} Let $A$ be a finite-dimensional pointed Hopf algebra with coradical 
$\ku \Gamma$, and let $R$ be the diagram of $A$, that is 
 $$\gr \mathcal A \simeq R \# \ku \Gamma,$$ 
and $R = \oplus_{n \ge 0} R(n)$ is a graded braided Hopf algebra in  $^{\Gamma}_{\Gamma}
\mathcal {YD}$ with $R(0) = \ku 1, R(1) = {\mathcal P}(R)$.

Assume that $R(1)$ is a Yetter-Drinfeld module of finite Cartan type with  braiding 
$(b_{ij})_{1 \leq i,j \leq \theta}$. For all i, let $q_{i} = b_{ii}, N_i = \text {ord}(q_{i})$. 
Assume that $\text{ord}(b_{ij})$ is odd and  $N_i$ is not divisible by 3 and $> 7$ 
for all $1 \leq i,j \leq \theta$.
\begin{enumerate}
\item For any $1 \leq i \leq \theta$ contained in a connected component of type $B_{n}$, 
$C_{n}$ or $F_{4}$ resp. $G_{2}$, assume that $N_{i}$ is not divisible by 5 resp. by 5 or 7. 
\label{degree2}
\item For any $1 \leq i,j \leq \theta$ and $i \nsim j$ assume $q_{i}q_{j} = 1$ 
or $ord(q_{i}q_{j}) = N_{i}$.\label{degree3}
\end{enumerate}

Then $R$ is generated as an algebra by $R(1)$, that is $A$ is generated 
by skew-primitive and group-like elements.
\end{teo}

\pf Let $S := R^*$  be the dual Hopf algebra of $R$ in the braided sense 
(see for example \cite[Section 2]{AG}). $S = \oplus_{n \ge 0} S(n)$ is a 
graded braided Hopf algebra in  $^{\Gamma}_{\Gamma}\mathcal {YD}$ with 
$S(0) = \ku 1, S(n)= R(n)^*$, for all $n \geq 0$. By assumption there are 
$h_{i} \in \Gamma, \eta_{i} \in \VGamma$, $ 1 \leq i,j \leq \theta$, with 
$b_{ij} = \eta_{j}(h_{i})$ for all $i,j$, and a basis $(y_{i})$ of $R(1)$ with 
$y_{i} \in R(1)_{h_{i}}^{\eta_{i}}$ for all i. Let $(x_i)$ in $V := S(1) = R(1)^*$ 
be the dual basis of $(y_i)$. Then $x_i \in V^{\chi_{i}}_{g_{i}}$ with $\chi_i = 
\eta_{i}^{-1} , g_i =h_{i}^{-1}$ and $b_{ij} = \chi_j(g_i) = \eta_{j}(h_{i})$ for 
all $1 \leq i,j \leq \theta$. Thus $V$ is a Yetter-Drinfeld module over $\Gamma$ 
with the same braiding as $R(1)$. By \cite[Lemma 5.5]{AS2}, $R$ is generated by  
$R(1)$ if and only if $S(1) = {\mathcal P}(S)$. Hence by duality,
$S$ is generated by $S(1)$, since $R(1) = {\mathcal P}(R)$. It is easy to see that 
$V = S(1) \subset 
{\mathcal P}(S)$. Hence there are canonical surjections of  graded braided Hopf algebras
$$T(V) \to S \to \toba(V).$$
Here $T(V)$ is the tensor algebra, the elements $x_i$ are primitive and of degree one, 
and both maps are the identity on $V$. The kernel $I$ of the first map is a homogeneous 
ideal generated by elements of degree $\geq 2$, a coideal and stable under the action 
and coaction of $\Gamma$. Since  $\toba(V) = T(V)/J$, where $J$ is the largest ideal 
with the same properties as $I$, there is a canonical surjection $S \to \toba(V)$.

The $x_i$ satisfy the Serre relations \eqref{serrebis} by Lemma \ref{Serre2}, and then 
the root vector relations \eqref{rootbis} by Lemma \ref{rootv}. Therefore it follows 
from the description of $\toba(V)$  in Theorem \ref{flk-tw} that $S = \toba(V)$. 
This means $S(1) = {\mathcal P}(S)$, hence by duality that $R$ is generated by $R(1)$.
\epf

A special case of the last Theorem together with a main result in \cite{AS2} allows 
to prove the following

\begin{cor}\label{degree1p} 
Let $p > 17$ be a prime number. Then any finite-dimensional pointed Hopf algebra
 with coradical $\ku (\Z/(p))^s$ for some natural number $s$ is generated by group-like 
 and skew-primitive elements.
\end{cor}

\pf Let $A$ be a finite-dimensional pointed Hopf algebra with coradical $\ku (\Z/(p))^s$ 
and let $R$ be the diagram of $A$. Then $R(1)$ is a Yetter-Drinfeld module of finite 
Cartan type by \cite[Corollary 1.2]{AS2}. Hence the claim follows from Theorem \ref{degree1}.
\epf

Let us state explicitly another Corollary of the Theorem.

\begin{cor} Under the hypothesis of Theorem \ref{degree1}, if the Dynkin diagram 
attached to the pointed Hopf algebra
is connected, then $A$ is generated by group-like and skew-primitive  elements. 
\qed \end{cor}

In principle, the idea behind the proof of Theorem  \ref{degree1} is as follows. Let 
$A$ be a finite-dimensional pointed Hopf algebra with coradical $\ku \Gamma$, 
$\Gamma$ any finite group. Let $R$ be the diagram of $A$, and $S := R^{*}$ the 
dual braided Hopf algebra. Consider the diagram $\widetilde{R}$ of the 
bosonization $S \# \ku \Gamma$. Then $\mathcal{P}(S)$ is naturally embedded in
 $\mathcal{P}(\widetilde{R})$ (and this embedding is in fact an isomorphism). 
 Moreover, $\text{dim}(\mathcal{P}(R)) \leq \text{dim}(\mathcal{P}(S)) \leq 
 \text{dim}(\mathcal{P}(\widetilde{R}))$, and
 $\text{dim}(\mathcal{P}(R)) = \text{dim}(\mathcal{P}(\widetilde{R}))$ if and 
 only if $S(1) = \mathcal{P}(S)$ or $R = \toba(\mathcal{P}(R))$.

Corollary \ref{degree1p} can also be seen as a direct consequence of Section 
\ref{serre} and \cite{AS2}: By \cite{AS2} $\mathcal{P}(\widetilde{R})$ is of 
finite Cartan type. Then the result follows from Theorem \ref{QSR} and 
\ref{powrootvec} applied to $A = S \# \ku \Gamma$.

The next theorem is another application of this principle. It shows that 
only very special dimensions are possible for finite-dimensional pointed Hopf algebras.

\begin{teo}\label{dimension}
For any finite group $\Gamma$ of odd order there is a natural number $n(\Gamma)$ such 
that the dimension of any finite-dimensional pointed Hopf algebra with 
coradical $\ku \Gamma$ is $\leq n(\Gamma)$.
\end{teo}

\pf Let $A$ be a finite-dimensional pointed Hopf algebra with coradical 
$\ku \Gamma$ and diagram $R$ and $\widetilde{R}$ as defined above. 
Since $R$ and $\widetilde{R}$ are braided Hopf algebras over 
$\Gamma$ of the same dimension, and $\text{dim}(\mathcal{P}(R)) 
\leq \text{dim}(\mathcal{P}(\widetilde{R}))$, we can iterate this 
process and after finitely many steps we obtain a graded braided
Hopf algebra $T$ over $\Gamma$ with $\text{dim}(R) = \text{dim}(T)$ and 
$T = \toba(\mathcal{P}(T))$.
By a result of Gra\~na \cite{Gn3} using \cite[Theorem 3.1]{AS2} which follows from 
\cite{L3}, the number of isomorphism classes of 
Yetter-Drinfeld modules $V$ over the fixed group $\Gamma$ with finite-dimensional 
$\toba(V)$ is finite. Thus we can take for $n(\Gamma)$ the product of the largest 
such dimension with the order of $\Gamma$.
\epf

\end{document}